\documentclass{amsproc}
\usepackage{amsmath}
\usepackage{amssymb}
\usepackage[matrix,poly,arc]{xy}

\usepackage{amsmath,amssymb,graphics}

\newtheorem{theorem}{Theorem}[section]

\theoremstyle{definition}
\newtheorem{definition}[theorem]{Definition}

\theoremstyle{remark}
\newtheorem{remark}[theorem]{Remark}
\newtheorem{remarks}[theorem]{Remarks}

\newtheorem{proposition}[theorem]{Proposition}
\newtheorem{conjecture}[theorem]{Conjecture}
\newtheorem{proposal}[theorem]{Proposal}
\newtheorem{observation}[theorem]{Observation}
\newtheorem{deftheorem}[theorem]{Definition/Theorem}

\numberwithin{equation}{section}

\def\gg {\mathfrak{g}}

\def\uu {\mathfrak{u}}

\renewcommand{\(}{\begin{equation}}
\renewcommand{\)}{\end{equation}}
\newcommand{\bea}{\begin{eqnarray}}
\newcommand{\eea}{\end{eqnarray}}
\newcommand{\R}{{\mathbb R}}
\newcommand{\C}{{\mathbb C}}
\newcommand{\Z}{{\mathbb Z}}
\newcommand{\Q}{{\mathbb Q}}

\def\oneone{\rlap 1\mkern4mu{\rm l}}

\def\S {\mathcal{S}}
\def\H {\mathcal H}
\def\A {\mathcal A}
\def\W {\mathcal W}
\def\V {\mathcal V}

\def\G {\mathcal G}

\def\L {\mathcal L}

\def\ee{\frak{e}}

\def\proof {{Proof.}\hspace{7pt}}

\def\endofproof {\hfill{$\Box$}\\}

\begin{document}

\title{Geometric and topological structures related to M-branes}

 \author{Hisham Sati}
  \address{Department of Mathematics, 
  Yale University, New Haven, Connecticut 06511\newline
Current address:
   Department of Mathematics, 
  University of Maryland, College Park, 
  Maryland 20742} 
  \email{hsati@math.umd.edu}
\subjclass[2010]{Primary 53C08, 55R65; 
Secondary  81T50,	55N20, 11F23}
\date{October 6, 2009
}

\keywords{String structures, Fivebrane structures,
$n$-bundles, differential cohomology,
K-theory, topological modular forms,  
generalized cohomology theories, 
anomalies, dualities}

\begin{abstract}
We consider the topological and geometric structures associated
with cohomological and homological objects in M-theory.
For the latter, we have  M2-branes and M5-branes, 
the analysis of which requires the underlying spacetime to 
admit a String structure and a Fivebrane structure, respectively. 
For the former, we study how the fields in M-theory are associated
with the above structures, with homotopy algebras, 
with twisted cohomology, and with 
generalized cohomology. We also explain how the corresponding 
charges should take values in Topological Modular forms.
We survey background material and related results in
 the process.




\end{abstract}

\maketitle

\tableofcontents

\section{Introduction and Setting}

String theory is concerned with a worldsheet $\Sigma^2$, usually a Riemann surface,
a target spacetime $M$, usually of dimension ten, and the space of maps 
$\phi: \Sigma^2 \to M$ between them. The study of the field theory on $\Sigma^2$ is the subject of
two-dimensional conformal field theory (CFT). The study of the maps is the 
sigma ($\sigma$-) model and the study of
the target space is the target theory where low energy limits, i.e. field and
supergravity theories, can be taken.

\vspace{3mm}
The target space theory involves fields called the Ramond-Ramond (RR) or the 
Neveu-Schwarz (NS) fields. These are differential forms or cohomology classes
which can be paired with homology cycles, the branes, namely D-branes 
and NS-branes respectively. These extended objects carry charges, generalizing
those carried by electrons.  
The RR fields are classified by K-theory \cite{MM} \cite{Wi1} and are twisted 
by the NS fields leading to a twisted K-theory classification \cite{Wi1} (see also 
\cite{BM}), in the sense of \cite{Ros} \cite{BCMMS}. Such a description has been
refined to more generalized cohomology theories, most notably elliptic cohomology. 
This was approached from cancellation of anomalies in type IIA string theory \cite{KS1}, 
studying the compatibility of generalized cohomology twists with S-duality in type IIB
string theory \cite{KS2}, and the study of modularity in the actions of type IIB string theory
and F-theory \cite{KS3} \cite{S4}.

\vspace{3mm}
The study of sigma models involves loop spaces as follows. 
The definition of spinors require lifting the classifying map of 
the tangent bundle from the special orthogonal group 
$\mathrm{SO}(n)$ to the spinor group 
$\mathrm{Spin}(n)$, which corresponds to killing the first
homotopy group $\pi_1(\mathrm{SO}(n))$ of $\mathrm{SO}(n)$. 
In string theory a further step is needed, 
namely lifting the spinor group to the 
String group $\mathrm{String}(n)$ by killing 
\footnote{This might perhaps more correctly be called ``cokilling" since
it corresponds to the Whitehead tower rather than to the 
Postnikov tower.} 
$\pi_3(SO(n))$, giving rise to String structures.
Existence of such structures is a condition for 
 the vanishing condition of the anomaly of a string
in the context of the index theory of Dirac operators on loop space \cite{Wit} \cite{Ki}.
The String structure is regarded as a lift of 
an $L \mathrm{Spin}(n)$-bundle
over the free loop space $LX$ through the Kac-Moody central extension 
$\widehat {L \mathrm{Spin}}(n)$-bundle
\cite{Ki} \cite{CP} \cite{PW} \cite{Mc}. 
This lift can also be interpreted 
as a lift of the 
original $\mathrm{Spin}(n)$-bundle down on target space $X$
to a principal bundle 
for the topological group
$\mathrm{String}(n)$ \cite{ST1}.  
This is, in fact, the realization of the nerve of a \emph{smooth} 
categorified group, 
the String 2-group \cite{BCSS} \cite{H}.
The above classification of String-bundles coincided with 
that of 2-bundles 
with structure 2-group the String 2-group \cite{BS}
\cite{BBK} \cite{BSt}. 
In addition to the above infinite-dimensional models, now
there is a finite-dimensional model for the String 2-group
\cite{SP}.
The elliptic genus is a loop space 
generalization of the $\widehat{A}$-genus as the index of the Dirac-Ramond operator
\cite{SWa} \cite{PSW} \cite{W} \cite{AKMW}. 
The Green-Schwarz \cite{GS} anomaly can be computed as essentially the 
elliptic genus \cite{LNSW}. Mathematically, the connection between elliptic 
genera and loop spaces has been studied, notably 
in \cite{And} and \cite{Liu}. The String structure , 
required by modularity, provides an orientation 
\cite{AHS} \cite{AHR} for 
TMF, the theory of Topological Modular Forms \cite{Ho} 
\cite{Go}.

\vspace{3mm}
At the level of conformal field theory, which is the quantum field theory 
on the worldsheet of the 
string, one has fields that are pulled back from the spacetime theory via 
the sigma model map, in addition to other fields. 
In two-dimensional supersymmetric quantum field theory the 
partition function, which is an integral modular 
function, is argued in \cite{ST1} \cite{ST2} to be an element in TMF. Another geometric 
description of elliptic cohomology via CFT is given in \cite{HK2}, which
builds on Segal's definition of CFT and on vertex operator algebras.  
A variant that has features of both \cite{ST1} and \cite{HK2} is proposed in 
\cite{GHK} using a newly-introduced notion of infinite loop spaces.


\vspace{3mm}
We emphasize two main points, central to the theme of this paper, 
concerning the above mathematical structures:
\begin{enumerate}
\item First, generally, that various structures appearing in this part of theoretical 
physics are much deeper (and thereby richer and more interesting) than the
sketchy physics literature about them indicates. 
\item Second, more specifically, that the above mathematical structures appearing 
in string theory are beginning to appear, even in a perhaps richer form, within the 
study of another theory, namely M-theory. 
\end{enumerate}
M-theory (cf.  \cite{Dynamics} \cite{Town})
is a conjectured theory in eleven dimensions that unifies all 
five ten-dimensional superstring theories. The theory is best understood through
these string theories and also via its classical low energy limit, eleven-dimensional
supergravity theory \cite{CJS}. Thus one strategy in studying the theory is to take 
eleven-dimensional supergravity and perform semi-classical quantization.
Due to quantum effects the process is only selectively reliable. Among the 
reliable terms are the topological terms, i.e. the terms that are not sensitive to 
the metric. Metric-dependent quantities might undergo drastic changes due to 
quantum gravitational effects. One way of keeping the metric requires 
taking some large volume limit, making sure that the scale is larger than 
the critical scale at which Riemannian geometry can no longer provide
a good description. As the theory is supersymmetrc, it will at least have a 
fermion (in this case, a section of the tensor product of the tangent bundle and the 
spin bundle), and since it involves gravity, it will also contain a metric, or graviton.
There is also an a priori metric-independent field, called the C-field. 
This is a higher-degree analog of a connection whose field strength -- the
analog of a curvature -- is denoted by $G_4$.

\vspace{3mm}
The fields (aside from the metric and fermions) 
in string theory and M-theory are differential forms at the rational 
level, i.e. at the level of description of the corresponding supergravity 
theories. Gauge invariance leads to a description of the fields in terms of 
de Rham cohomology. 
Quantum mechanically, these generically become 
integral-valued and hence one needs to go beyond de Rham 
cohomology to integral cohomology. Dually, these fields can be 
described via homology cycles that admit extra geometric structures,
such as Spin structures and vector bundles. This dual homological 
picture is captured by the notion of {\it branes}, namely D-branes in 
string theory and M-branes in M-theory. The fields and the branes 
are of specified dimensions, determined by the corresponding 
theory. In particular, D-branes have odd (even) spacetime dimensions,
and hence even (odd) spatial dimensions, for type IIA (type IIB) 
string theory. For M-branes, spatial dimensions two for the M2-brane
 \cite{BST} and five for the M5-brane \cite{Gu} occur.

\vspace{3mm}
The five string theories in ten dimensions are related through a web of dualities
(see e.g. \cite{Lec} for a survey).
The first kind of duality is called T-duality ( ``T" for Target space), 
which relates two different theories
on torus bundles, where the first theory with fiber a torus is related to a second 
theory with a fiber the dual torus. An example of this is T-duality between type IIA
and type IIB string theories.
The second kind of dualities is S-duality,
or strong-weak coupling duality, which relates a theory at a high value of 
some parameter (strong coupling) to another theory at a low value of the same
parameter (weak coupling).
This generalizes the usual electromagnetic duality
between electric and magnetic fields in four dimensions. Whenever we discuss
dualities in this paper we will focus mostly on S-duality. 
A duality within the same theory is called self-duality. An example of this is the
self-duality among the RR fields in type IIA/B string theory. At the rational 
level this is simply a manifestation of Hodge duality. A delicate discussion of 
such matters can be found in  \cite{Fr} \cite{FMS}.
Another important example of this is S-(self)duality in type IIB 
string theory. Some subtleties on the relation of this duality to generalized 
cohomology are discussed in \cite{BEJMS} \cite{KS2}.

\vspace{3mm}
The above dualities can be most directly seen at the level of fields. 
By Poincar\'e duality, such dualities also manifest themselves at the
level of branes. Branes of even (odd) spatial dimensions in type IIA (type IIB)
 are dual to odd (even) branes in 
type IIB (IIA) string theory. Furthermore, within the same type II theory,
D-$p$-branes are dual to D-$(6-p)$-branes, $p$ even or odd for type 
IIA or type IIB, respectively. This is a homological manifestation of the 
self-duality of the RR fields. There is a similar duality in eleven dimensions
which relates  the $C$-field $C_3$ to its dual $C_6$, which, at the rational level, 
is Hodge duality between the corresponding fields strengths $G_4=*G_7$. 
This duality on the fields also has a homological interpretation as a duality
between the M2-brane and the M5-brane.

\vspace{3mm}
String theories in ten dimensions can be obtained
from eleven-dimensional M-theory via dimensional reduction and/or duality transformations.
M-theory on the total space of a circle bundle gives rise to type IIA string 
theory on the base space.
By pulling back the fields along the section of the 
circle bundle $\pi$, assumed trivial, one gets 
a D2-brane from an M2-brane \cite{Town2}
and a NS5-brane from an M5-brane.
On the other hand, upon integration over the fiber of
$\pi$, the
M2-branes give rise to strings \cite{DHIS} 
and the M5-brane give rise to D4-branes 
\cite{Town2}.  
The branes of type IIB string theory can also be 
obtained from those of M-theory on a torus bundle 
\cite{Sch}.
Similar relations hold at the level of fields: Integrating 
$G_4$ in eleven dimensions over the circle gives 
$H_3=\pi_*(G_4)$, the field strength of the NS $B$-field. 
On the other hand, pulling back $G_4$ along
a section $s$, again assuming the circle bundle is trivial,
gives a degree four field $F_4=s^*G_4$, which is one
component of the total RR field. An invariant
description of 
this is given in 
 \cite{FS} and \cite{TDMW}.

\vspace{3mm}
Branes carry charges-- a notion that can be made mathematically precise--
that can be viewed either as classes of bundles in generalized cohomology
or as their images in rational or integral cohomology under a (normalized) 
version of the Chern character map.   
A working mathematical definition of D-branes and their charges can be
found in \cite{BMRS}:
  A D-brane in ten-dimensional spacetime $X$ is a triple 
  $(W,E, \iota)$, where $\iota : W\hookrightarrow X$ is a closed, embedded 
  submanifold and $E \in {\rm Vect}(W)$ is a complex vector bundle over $W$. 
  The submanifold $W$ is called the worldvolume and $E$ the Chan-Paton bundle of the
  D-brane.
The charges of the D-branes \cite{Pol}
can be classified, in the absence of the NS fields, by
K-theory of spacetime \cite{MM}, namely by
$K^0(X)$ for type IIB \cite{Wi1} and by $K^1(X)$ for type IIA \cite{Ho}.
The fields are also classified by K-theory of spacetime but 
with the roles of $K^0$ and $K^1$ interchanged \cite{MW} \cite{FH}. 
In the presence of the NS B-field, 
or its
field strength $H_3$, the relevant K-theory is 
twisted K-theory, as was shown in \cite{Wi1} \cite{FW} \cite{Kap} by analysis of
worldsheet anomalies for the case the NS field  
$[H_3] \in H^3(X, \mathbb Z)$ is a torsion class, and in \cite{BM} for the
nontorsion case. 
Twisted K-theory has been studied for some time \cite{dk}  \cite{Ros}.  
More geometric flavors were given in \cite{BCMMS}. Recently, the theory
was fully developed by Atiyah and Segal \cite{AS1} \cite{AS2}. 
The identification of twisted D-brane charges with elements in twisted
K-theory requires a push-forward map and a Thom isomorphism in the latter,
both of which are established in \cite{CW}.

\vspace{3mm}
The study of the M-brane charges and their relation to (generalized) cohomology
is one of the 
main goals of this paper. 
Given the relation between generalized cohomology and string theory on one hand,
and between string theory and M-theory on the other, it is natural to ask whether elliptic cohomology or 
$TMF$ have to do with M-theory
directly. Suggestions along these lines have been given in \cite{KS1} \cite{S1} \cite{S2}
\cite{S3}. In particular in \cite{KS1} it was proposed that the M-branes should be
described by TMF in the sense that the elliptic refinement of the partition function
originates from interactions of M2-branes and M5-branes. Furthermore,
in \cite{S3} it was observed that the M-theory field strength $G_4$, rationally,
 can be viewed as part of a twist of the de Rham complex and suggested 
 that the lift to generalized cohomology would be related to a twisted version 
 $TMF$. A twisted differential is one of the form $d + \alpha\wedge$ 
acting on differential forms, where $d$ is the de Rham differential
and $\alpha$ is a differential form. Twisted rational cohomology is then 
the kernel  modulo the  image of $d + \alpha\wedge$. 
When $\alpha$ is a 3-form then one gets the image of twisted K-theory 
under the twisted Chern character (cf. \cite{AS2} \cite{BCMMS} \cite{MSt}).

\vspace{3mm}
The field strength $G_4$ on an eleven-manifold $Y^{11}$
 satisfies the shifted quantization condition \cite{Flux}
\(
[G_4] +  \frac{1}{2}\lambda(Y^{11}) = a \in H^{4}(Y^{11}; \Z)\;,
\label{shift}
\)
with $\lambda(Y^{11})=\frac{1}{2}p_1(Y^{11})$, where $p_1(Y^{11})$ is the first Pontrjagin class of the tangent bundle
 $TY^{11}$ of $Y^{11}$ and
$a$ is the degree four class that characterizes an $E_8$ bundle in
M-theory \cite{DMW}. 
There is a one-to-one
correspondence between $H^4(M, \Z)$ and isomorphism classes of principal $E_8$
bundles on $M$, when the dimension of $M$ is less than or equal to 15, which is 
the case for $Y^{11}$ in M-theory. 
This follows from homotopy type  of $E_8$ 
being of the form $(3, 15, \cdots)$.
Up to the 14th-skeleton 
$E_8$ is homotopy equivalent to the 
Eilenberg-MacLane space $K(\Z,3)$ so that
up to the 15th-skeleton 
the classifying space $BE_8$ is $\sim K(\Z,4)$.
For the homotopy classes of
maps 
$[M,E_8]= [M, K(\Z,3)] = H^3(M,\Z)$  
if $\mbox{dim}M \leq 14$, and similarly $\{$Equivalence classes of $E_8$
bundles on $M \}=[M,BE_8]= [M, K(\Z,4)]=
H^4(M, \Z)$ if $\mbox{dim}M \leq 15$. 
Therefore, corresponding to an element $a \in H^4(M,\Z)$ we have an $E_8$
principal bundle $P(a)\to M$ with $p_1(P(a))= a$.

\vspace{3mm}
In \cite{Wa} the notion of a twisted String structure was defined, where 
the twist is given by a degree four cocycle. This degree four generalization of 
the twisted ${\rm Spin}^c$ structure in degree three
was anticipated there to be related to the flux quantization condition (\ref{shift}).
This was made explicit and precise in \cite{SSS3}, where also the Green-Schwarz anomaly
was shown to be more precisely the obstruction to having a 
refined twisted String structure.

\vspace{3mm}
The study the dual NS- and M-theory fields of degree eight leads to even higher structures. 
The first appears in a dual form of the Green-Schwarz 
anomaly cancellation. The second satisfies a condition analogous to
(\ref{shift}) in degree eight, as shown in \cite{DFM}. 
In \cite{SSS1} Fivebrane structures were introduced and the corresponding 
differential geometric structures, i.e. the higher bundles, are constructed. 
They are systematically studied in \cite{SSS2} and shown to emerge 
within the description of the 
dual of the Green-Schwarz anomaly, involving the Hodge dual of $H_3$ 
in ten dimensions, as well as in the description of the dual of the M-theory 
 $C$-field in eleven dimensions. 
Fivebrane structures can be twisted in the same way that 
String structures can. Twisted Fivebrane structures are defined and studied 
in \cite{SSS3}, and their obstructions can be matched with the 
dual of the Green-Schwarz anomaly. In fact, both twisted
String and Fivebrane structures are refined in \cite{SSS3}
to the differential case, using a generalization of the discussion in
 \cite{HS}. The notion of (twisted) String and, to some extent, 
Fivebrane structures can in fact be described
in various ways, via:

\noindent {\bf 1.}  Principal and associated bundles.

\noindent {\bf 2.} Gerbes and differential characters. 

\noindent {\bf 3.} {\v C}ech cohomology and Deligne cohomology.

\noindent {\bf 4.} Loop bundles. 

\noindent {\bf 5.} 2-bundles and 6-bundles and their  
2-algebras and 6-algebras, respectively.

\vspace{3mm}
Another purpose of this paper is to provide the generalized cohomology aspect.
The appearance of $\lambda$ in (\ref{shift}) and the subsequent 
interpretation in terms of twisted String structures suggests a relation 
to a theory that admits that structure as an orientation. 
A Spin manifold $M$ has a characteristic class $\lambda$ such that
$2\lambda = p_{1} (M).$  The paper \cite{AHS} shows that $M$ admits a $TMF$ orientation
if $\lambda = 0$.  More precisely, a String structure on a Spin
manifold $M$ is a choice of trivialization of $\lambda$, and in \cite{AHS}
 it is shown
that a String structure determines an orientation of $M$ in
$TMF$-cohomology. 
In this paper we argue that the shifted quantization condition (\ref{shift}) provides
a twist for $TMF$ and hence that (\ref{shift}) defines a twist for $TMF$.
More precisely, the cohomology class of $G_4$ is $\frac{1}{2}\lambda -a$,
which we view as a twist of the $TMF$ orientation by the degree four class 
$a$ of the $E_8$ bundle. Since $E_8 \sim K(\Z, 4)$ up to dimension 
14 then $a$ a priori can be any class in $K(\Z, 4)$. However, fixing an
$E_8$ bundle completely fixes $a$. Conjectures for using a twisted form
of $TMF$ and elliptic cohomology -- and hence that such structures should exist--
to describe the fields in M-theory go back to \cite{S2} \cite{S3}.

\vspace{3mm}
Given the interpretation of $G_4$ and its dual in terms of twisted String 
and Fivebrane structures 
and the proposed connection to twisted $TMF$, it is natural to consider
the corresponding homological objects.  The charge of the M5-brane is,
rationally, the 
value of the integral of $G_4$ over
the unit sphere in the normal bundle in the eleven-dimensional manifold.
We interpret the charge of the M5-brane in full and not just rationally, 
in a Riemann-Roch setting, as the direct image of 
an element in twisted $TMF$ on the M-brane. 
The interpretation of the charges as
elements of twisted $TMF$ 
uses some recent work \cite{Atwist} on
push-forward and Thom isomorphism in $TMF$.
This is a higher degree generalization of the case of D-branes, 
where
the $H$-field defines a twist for the ${\rm Spin}^c$ structure
and the charge is defined using the push-forward and Thom
isomorphism in twisted K-theory.

\vspace{3mm}
Given the topological structures defined by $G_4$ and its dual, it is natural
to ask for geometric models for the corresponding potentials, i.e. the $C$-field $C_3$ 
and its dual. There is the $E_8$ model of the $C$-field \cite{DMW} \cite{DFM},
mentioned above,
which is essentially a Chern-Simons form, or more precisely a shifted
differential character, where the shift on the $E_8$ bundle class 
$a$ is given by the factor 
$\frac{1}{2}\lambda$. In an alternative model in terms of 2-gerbes
\cite{AJ}, picking a
connection $A$ on the $E_8$ principal bundle $P(a)$ gives
 a Deligne class, the Chern-Simons
2-gerbe  $CS(a)$, with $a$ its characteristic class.
The interpretation we advocate in terms of twisted 
String structures allows for other interpretations of the $C$-field, 
where the term $\frac{1}{2}\lambda$ is now the `main' term, since
it is responsible for defining the String structure, and the term $a$ 
coming from the $E_8$ bundle is merely a twist for that structure. 
This gives the tangent bundle-related term a more prominent role.
In fact, this model can be described at the differential twisted 
cohomology level \cite{SSS3}. In this paper we provide an identification 
of the $C$-field as the {\it String class}, i.e. as the class that provides the 
trivialization of the String structure.  We also do the same for the dual 
of the $C$-field which we identify with the {\it Fivebrane class}, 
the class that provides the trivialization for the Fivebrane structure. 
We also give an alternative description using differential 
characters.

\vspace{3mm}
The appearance of higher chromatic phenomena in \cite{KS1} \cite{KS2}
\cite{KS3} in relation to string theory, and the appearance of 
Fivebrane structures in degree eight in string theory and M-theory 
naturally leads to the question of whether higher degree twists 
exist in this context. Indeed, in \cite{7twist} it was shown that a degree 
seven twist occurs in heterotic string theory, manifested via a differential 
of the form $d + H_7 \wedge$. The lift of this twisted rational cohomology
to generalized cohomology suggests the appearance of the second 
generator $v_2$ at the prime 2 in theories descending from complex cobordism
$MU$. This generalizes the situation in degree three, where the Bott generator
$u=v_1$ in K-theory appears via $d + v_1 H_3\wedge$. 
In this paper we consider the fields in M-theory as part of a twist in
de Rham cohomology, extending and refining 
the discussion in \cite{S1} \cite{S2} \cite{S3}. In addition we consider 
duality-symmetric twists, i.e. twisted differentials whose twists are 
uniform degree combinations of the $H$-field and its dual
as well as $G_4$ and its dual. The second case leads to an interesting 
appearance of the M-theory gauge algebra, which in turn leads to 
the super-tranlation algebra. We also provide an $L_{\infty}$-algebra
 description of this
gauge algebra.

\vspace{3mm}



\begin{table}[h]
  \begin{center}
  \begin{tabular}{lr|ll}
    \multicolumn{2}{c}{\bf String Theory}
    &
    \multicolumn{2}{c}{ \bf M-Theory}
    \\
    \hline
    &
    \\
    sigma model  $\phi : \Sigma^2 \hookrightarrow X^{10}$ &
    & sigma model  $\Phi : {M}^3 \hookrightarrow Y^{11}$ &
 \\
    \hline
    &
    \\
 $\psi\in \Gamma( S\Sigma^2 \otimes \phi^* TX^{10})$, 
$B_2$ 1-gerbe & &
$\psi \in \Gamma( S{M}^3 \otimes {\mathcal N}({M}^3 \hookrightarrow Y^{11})$, 
$C_3$ 2-gerbe 
 \\
    \hline
    &
    \\
D-brane $\supset$ $\partial \Sigma^2$ 
& & M5-brane  $\supset$ $\partial {M}^3$ &
 \\
    \hline
    &
    \\
Freed-Witten condition $W_3 + [H_3]=0$ &
& Witten Flux quantization $\frac{1}{2}\lambda + [G_4] = a$ &
\\
\hline
\end{tabular}
  \end{center}
  \caption{
   Extended objects in string theory and in M-theory.   
  }
\end{table}

\vspace{-3mm}
This paper is written in an expository style, even though it is mainly 
about original research. In fact, there are three types of material:
\begin{itemize}
\item Survey of known results, with some new perspectives a well as
providing some generalization. 
\item New research established here.
\item New research announced and outlined 
here and to be more fully developed in the future. This material
is mostly based on discussions with Matthew Ando, Chris Douglas, 
Corbett Redden, Jim Stasheff, and Urs Schreiber. 
\end{itemize}
 We hope that the expository style makes it more self-contained 
and accessible to mathematicians interested in this area of interaction 
between physics on one hand and geometry and topology on the other.

\section{The M-theory $C$-Field and String Structures}
Consider an eleven-dimensional Spin manifold $Y^{11}$ with metric $g$.
Corresponding to the tangent bundle $TY^{11}$ with structure group $SO(11)$
 there is the Spin bundle $SY^{11}$ with structure group the Spin 
 group ${\rm Spin}(11)$, the double cover of $SO(11)$. Let $\omega$ be 
 the Spin connection on $SY^{11}$ with curvature $R$. Using the metric we 
 can identify the cotangent bundle $T^* Y^{11}$ with the tangent bundle.
 The field content of the
 theory is the metric $g$, a spinor one-form  $\psi$ , i.e. a section 
 $\in \Gamma( SY^{11} \otimes T^*Y^{11})$, and a degree three form $C_3$.

\subsection{The Quantization Condition and the $E_8$ model for the $C$-field}

The quantization condition on $G_4$ on $Y^{11}$ is given in equation
(\ref{shift}). 
This can be obtained either from the
partition function of the membrane or from the reduction to the 
heterotic $E_8 \times E_8$ theory on the boundary \cite{Flux}. 
We will reproduce this result using the membrane partition function
in a fashion that is essentially the same  
as appears in \cite{Flux} but with details retained.  
Since the Spin cobordism groups $\Omega^{\rm Spin}_{4k+3}$ are zero, 
we can extend both the membrane worldvolume $M^3$ and the target
spacetime $Y^{11}$ as Spin manifolds to bounding manifolds
$X^4$ and $Z^{12}$, respectively. In fact we can also extend $E_8$
bundles on $M^3$ and $Y^{11}$ to $E_8$ bundles on $X^4$ and 
$Z^{12}$, respectively, since $M{\rm Spin}_{i}\left( K(\Z, 4) \right)=0$
for $i=3, 11$ \cite{Stg} where $BE_8$ has the homotopy of $K(\Z, 4)$ in our range of 
dimensions. The effective action involves two factors:
\begin{enumerate}

\item The `topological term' $\exp i \int_{M^3} C_3$, 

\item The fermion term $\exp i \int_{M^3} \overline{\psi} D \psi$, where the integrand in the
exponential is the pairing in spinor space of the spinor $\psi$ with the spinor 
$D\psi$, where $D$ is the Dirac operator.

\end{enumerate}
The first factor will simply give $\exp i \int_{X^4} G_4$. 
Now consider the second factor. Corresponding to the map
$\phi: X^4 \to Z^{12}$ we have the index of the Dirac operator for
spinors that are sections of $S(X^4) \otimes \phi^* TZ^{12}$ given
via the index theorem by the degree two expression
\bea
{\rm Index} D &=& \left[ \int_{X^4} \widehat{A}(X^4) \wedge {\rm ch}(\phi^*TZ^{12})  
\right]_{(2)}
\nonumber\\
&=& \int_{X^4} \left[ 1 - \frac{1}{24} p_1(TX^4) \right] 
\left[ {\rm rank}(TZ^{12}) + {\rm ch}_2(\phi^*TZ^{12}) \right]
\nonumber\\
&=& \int_{X^4} \frac{1}{2} p_1(\phi^*TZ^{12}) - \frac{1}{2} p_1(X^4)\;.
\eea 
On the other hand we have a split of the restriction of the
tangent bundle $TZ^{12}$ to $X^4$ as
$
TZ^{12}=TX^4 \oplus {\mathcal N}X^4
$
with ${\mathcal N}X^4$ the normal bundle of $X^4$ in $Z^{12}$.
Taking the characteristic class $\lambda$ of both sides  
 we get that the index is equal to $\lambda({\mathcal N}X^4)$.
The effective action involves the square root of the index so that
the contribution from the second factor in the effective action is
$
\exp 2 \pi i \frac{1}{2} \lambda({\mathcal N}X^4)
$.
This gives the result in \cite{Flux} provided we assume that
$
{\mathcal N}( M^3 \hookrightarrow Y^{11}) \cong 
{\mathcal N}(X^4 \hookrightarrow Z^{12})
$.

\begin{proposition}
The M-theory field $G_4$ satisfies the quantization condition
(\ref{shift}).
\end{proposition}
\begin{remarks}
The  condition (\ref{shift}) has the following consequences \cite{Flux}:
\begin{enumerate}
\item When $\lambda$ is divisible by two then $G_4$ cannot be set to 
zero.

\item When $\frac{\lambda}{2}$ is integral then the tadpole anomaly of
\cite{SVW} vanishes.

\item The relation is parity-invariant, i.e. requiring $G_4 - \frac{\lambda}{2}$
to be integral is equivalent to $G_4 + \frac{\lambda}{2}$ being integral.
\end{enumerate}
\end{remarks}

  \vspace{3mm}
\noindent {\bf The $E_8$ model of $C$-field.}
The $C$-field at the level of supergravity will be just a real-valued 
three-form $C_3 \in \Omega^3 (Y^{11}; \R)$. The field strength is
$G_4= dC_3 \in \Omega^4(Y^{11})$. This is invariant under gauge 
transformations $C_3 \to C_3 + d\phi_2$, where $\phi_2$ is a
two-form. Factoring out by the gauge transformations amounts
to declaring the fields to be in cohomology. Upon quantization, several 
features become important: integrality via holonomy, the 
presence of torsion, and possible appearance of anomalies. Taking
these into account, a model for the $C$-field was obtained in 
\cite{Flux} and further developed in 
\cite{DMW} \cite{DFM} \cite{Gauss}. 

\vspace{3mm}
Let $P$ be a principal $E_8$ bundle over $Y^{11}$ with the characteristic
class $a$ pulled back from $H^4(BE_8; \Z)$. Let $A$ be a
connection on $P$ with curvature two-form $F$. The $C$-field in this model
is given by 
$
C= CS_3(A) -\frac{1}{2} CS(\omega) + c
$,
where $CS_3(A)$ is the Chern-Simons invariant for the connection $A$,
$CS(\omega)$ is the Chern-Simons invariant of the connection $\omega$
on the Spin bundle,
and $c$ is the harmonic representative of the $C$-field which dominates
at long distance approximation, i.e. in the supergravity regime. The 
Chern-Simons forms and the Pontrjagin forms are related as
$
dCS_3(A)= {\rm Tr} F\wedge F
$,
$dCS_3(\omega) = {\rm Tr} R \wedge R
$,
so that the field strength of the $C$-field is given by
\(
G_4= {\rm Tr} F\wedge F - \frac{1}{2} {\rm Tr} R \wedge R + dc\;. 
\label{g4 formula}
\)
The cohomology classes are 
$
\left[ {\rm Tr} F \wedge F\right]_{DR} = a_{\R}
$,
$\left[ {\rm Tr} R \wedge R\right]_{DR} = \frac{1}{2} \left( 
p_1(TY^{11})\right)_{\R}
$.
Globally, the $C$-field can be described as the pair \cite{DFM}
\(
(A, c) \in {\mathcal E}_P (Y^{11}) := {\mathcal A} (P(a)) \times \Omega^3(Y^{11}),
\)
where ${\mathcal A}(P(a))$ is the space of smooth connections on 
the bundle $P$ with class $a$.  


\subsection{Twisted String structures}

 \begin{definition}
 An $n$-dimensional manifold $X$ admits a
 \emph{ String structure} 
 if the classifying
 map $X \to BO(n)$ of the tangent bundle $T X$ lifts to the
 classifying space
 $B{\rm String} :=B\mathrm{O}\langle 8 \rangle$.
\(
  \raisebox{20pt}{
  \xymatrix{
    &
     \mathrm{B O} (n) \langle 8 \rangle 
     \ar[d]
     \\
     M
     \ar[r]_f
     \ar@{..>}[ur]^{\hat f}
     &
     \mathrm{B O}(n)~.
  }
  }
\)
\end{definition}

\begin{remarks}
\begin{enumerate}
\item The obstruction to lifting a Spin structure on $X$
    to a String structure on $X$ is the fractional
    first Pontrjagin class
    $ \frac{1}{2}p_1(T X)$.
    \item The set of lifts, i.e. the set of String 
structures for a fixed Spin structure
is a 
torsor for a quotient  
of the third integral cohomology group $H^3(X; \Z)$. 
\end{enumerate}
\end{remarks}

\begin{definition}[\cite{Wa} \cite{SSS3}]
  An $\alpha$-twisted String structure
  on a brane $\iota : M \to X$
  with Spin structure classifying map $f : M \to BO(n) \langle 4 \rangle$ is a
  cocycle $\alpha : X \to K(\mathbb{Z}, 4)$ and a map $c: B O (n) \langle  4 \rangle 
  \to K(\Z, 4)$ such that there is a homotopy $\eta$
  between  $B O (n) \langle  4 \rangle $ and $X$, 
  $$  
  \raisebox{20pt}{
    \xymatrix{
       M
       \ar[rr]^f_>{\ }="s"
       \ar[d]_\iota
       &&
       B O (n) \langle  4 \rangle
       \ar[d]^{c}
       \\
       X
      \ar[rr]_\alpha^<{\ }="t"
       &&
       K(\mathbb{Z}, 4)
            \ar@{=>}^\eta "s"; "t"
   }
    }\;.
    $$
The homotopy $\eta$ is the coboundary that relates the cocycle 
$cf$ to the cocycle $\alpha \iota$. Hence on cohomology classes it says that
the fractional Pontrjagin class of $M$ does not necessarily vanish, but is equal to 
the class $\iota^* [\alpha]$. 
\end{definition}

\vspace{3mm}
The twisted String case was originally considered in \cite{Wa}.
The definition is refined to structures dubbed 
$\mathcal{F}^{\langle m \rangle}$ 
that account for obstructions that are fractions of the ones for the 
String structures. For example, the fractional class
$\frac{1}{4}p_1$  shows up in (\ref{shift}). 
As an application, which was originally the motivation:
\begin{theorem} [{\cite{SSS3}}]
\begin{enumerate}
\item The Green-Schwarz anomaly cancellation condition defines a twisted String structure 
pulled back from $BO(10) \langle 4 \rangle =B {\rm Spin}(10)$. The twist $\alpha$ in this case 
is given by (minus) the degree four class of the $E_8 \times E_8$ bundle.
 \\
\item The anomaly cancellation condition in heterotic M-theory and
the flux quantization condition in M-theory each define a
twisted String structure pulled back from
 $\mathcal{F}^{\langle 4\rangle}=BO\langle \frac{1}{4}p_1 \rangle$.
The twist $\alpha$ in this case is given by $[G_4]$ minus 
the class of the $E_8$ bundle.
\end{enumerate}
\end{theorem}
The division of $\lambda$ by 2 require some refinement of the structure
\cite{SSS3} as mentioned in Remarks \ref{rem 4}.
The relation to orientation in generalized cohomology is 
discussed in the study of the membrane partition function in
section \ref{st mem}. Note that when we consider M-theory with
a boundary $\partial Y^{11}$, where essentially 
the heterotic string theory is defined, 
$[G_4]$ would be zero when restricted to $\partial Y^{11}$. In this
case the flux quantization condition defines a $a$-twisted String 
structure $\frac{1}{2}p_1(\partial Y^{11}) =a |_{\partial Y^{11}} \neq 0$
on that boundary. This is discussed further in section 
\ref{M ch tmf}.

\subsection{The $C$-field as a String class}
\label{c as string}

Let $\langle \cdot , \cdot \rangle$ be a suitably normalized Ad-invariant metric 
on the Lie algebra $\frak{spin}(11)$ of ${\rm Spin}(11)$. Then the 4-dimensional
Chern-Weil form $\langle R \wedge R \rangle \in \Omega^4(Y^{11})$
on $S(Y^{11})$ is one-half the first Pontrjagin class (restoring normalization)
$
\frac{1}{2}p_1(S, \omega) =- \frac{1}{16 \pi^2} {\rm Tr} (R \wedge R)
$.
We will later assume the following condition on the Pontrjagin class 
$
\frac{1}{2}p_1(Y^{11})=0$ in $H^4(Y^{11};\R)
$.
This means that the bundle $S(Y^{11})$ admits a String structure. A choice of String structure
is given by a particular cohomology class 
$\S \in H^3\left( S(Y^{11});\Z \right)$. This element restricts to the fiber as the standard 
generator of $H^3({\rm Spin}(11); \Z) \cong \Z$. 

\vspace{3mm}
In terms of the curvature $R$ of the connection
$\omega$ on $S(Y^{11})$, the condition $\frac{1}{2}p_1(S(Y^{11}))=0 \in H^4(Y^{11};\R)$
means that 
$
\left[ {\rm Tr} R \wedge R \right] =0 \in H^4(Y^{11};\R)
$.
The Chern-Simons 3-form of the connection $\omega$ on the principal
bundle $S(Y^{11})$ is the right-invariant form
$
CS_3(\omega) := \langle \omega \wedge R \rangle 
-\frac{1}{6} \langle \omega \wedge [\omega, \omega] \rangle
\in \Omega^3(S)
$, 
whose pull-back to the fiber via the inclusion map $i: {\rm Spin}(11) \hookrightarrow
S$ is the real cohomology class
$
-\frac{1}{6} \langle \omega \wedge [\omega, \omega] \rangle \in H^3({\rm Spin}(11);\R)
$
associated to the real image of the standard generator of $H^3({\rm Spin}(11);\Z) \cong \Z$.
The Chern-Simons form $CS_3(\omega)$ provides a trivialization for the
zero cohomology class above:
$
dCS_3(\omega) = {\rm Tr} R \wedge R \in \Omega^4(S(Y^{11}))
$.

\vspace{3mm}
Now we consider the geometry on the total space of the Spin bundle. Corresponding 
to a choice of Riemannian metric $g$ on $Y^{11}$ and a connection 
$\omega$ on $S(Y^{11})$ gives rise to a metric $g_S$
on the 
total space $S(Y^{11})$. Under the decomposition 
$TS(Y^{11}) \cong \pi^* (TY^{11} \oplus {\frak{spin}}(11))$ of the tangent space
into orthogonal vertical and horizontal subspaces, the metric decomposes
as
$
g_S : = \pi^* (g \oplus g_{{}_{{\rm Spin}(11)}} )
$,
where $g_{{}_{{\rm Spin}(11)}}$ is the metric on the fiber.

\vspace{3mm}
In relating the fields on the base to classes on the 
total space, one is forced to use the adiabatic limit, introduced in \cite{CMP}, 
of the metric on the 
total space. In particular, as is constructed in \cite{Red}, 
the String structure $\S$ on $S(Y^{11})$ is related
to a form on the base in this fashion. This way there is a one-parameter family of 
metrics 
$
g_{\delta} = \pi^* \left( \frac{1}{\delta^2} g \oplus g_{{}_{{\rm Spin}(11)}}
\right)
$
on the bundle $S(Y^{11})$ with parameter $\delta$.
This is reminiscent of a Kaluza-Klein ansatz
frequently used in supergravity. Now consider the adiabatic limit 
$\delta \to 0$ of $g_{\delta}$.
Note that metrics in the adiabatic limit have
been used in this form e.g. in \cite{TDMW} in the reduction from M-theory
in eleven dimensions to type IIA string theory in ten dimensions.

\vspace{3mm}
\noindent{\bf Relative Chern-Simons form on $Y^{11}$.}
Before considering the $C$-field we will need the following.
Consider an $E_8$ bundle $P$ over
$
Y^{11}
$
with connection $A \in \Omega^1(P; \frak{e}_8)$ and curvature 
$
\pi^* \Omega = dA + \frac{1}{2} [A, A] \in \Omega^2(P; \frak{e}_8),
$
where $\Omega \in \Omega^2(Y^{11}; {\rm ad}P)$, a two-form with
values in ${\rm ad}P$. 
The space $\A_P$ of connections on $P$ is an affine space modeled on
$\Omega^1(Y^{11}; {\rm ad}P)$. Using \cite{CSII} we can define a relative
Chern-Simons invariant on the base. Given two connections $A_0$ and
$A_1$ in $\A_P$, the straight line path 
$A_t=(1-t) A_0 + tA_1$, $0 \leq t \leq 1$, 
determines a connection ${\mathbb A}$ on the bundle
$
\xymatrix{
[0,1] \times P 
\ar[r]
&
[0,1] \times Y^{11}
}
$.
The relative Chern-Simons form is then
\(
CS_3(A_1, A_0):=- \int_{[0,1]} {\rm Tr} F^2({\mathbb A}) \in \Omega^3(Y^{11})
\label{rel CS}
\)
Using Stokes' theorem
\bea
dCS_3(A_1, A_0) &=& -d \int_{[0,1]} {\rm Tr} F^2 ({\mathbb A})
\nonumber \\
&=& - \int_{[0,1]} d {\rm Tr} F^2 ({\mathbb A}) + (-1)^{(11-4)}
\int_{\partial [0,1]} {\rm Tr} F^2({\mathbb A})
\nonumber \\
&=&
0 + {\rm Tr} F^2 (A_1) - {\rm Tr} F^2(A_0)\;.
\eea
Now we are ready to consider the $C$-field.

\vspace{3mm}
\noindent {\bf Invariance of the $C$-field.}
The $C$-field is invariant under the following transformations \cite{DFM}:
$A'= A + \alpha$ and
$C_3'=C_3 - CS_3(A, A + \alpha) + \Lambda_3$,
where $\alpha \in \Omega^1\left( {\rm ad} S(Y^{11}) \right)$
and $\Lambda_3$ is a closed 3-form on $Y^{11}$. The Chern-Simons invariant 
takes values in $\R/\Z$ so that it is not defined as a differential form unless 
exponentiated. The relative Chern-Simons invariant is defined as in 
(\ref{rel CS}).

\vspace{3mm}
If we include the Spin bundle $S(Y^{11})$ then we also have an invariance
of the connection $\omega$ and a corresponding shift in the relative Chern-Simons
form of $\omega$. Thus we have

\begin{proposition}
The $C$-field is invariant under the following transformations
\begin{enumerate}
\item $\omega'=\omega + \beta$,
\item $A' = A + \alpha$,
\item $C_3'= C_3 - CS_3(A, A+ \alpha) + \frac{1}{2} CS_3(\omega, \omega + \beta)
+ \Lambda_3$,
\end{enumerate}
where $\alpha \in \Omega^1({\rm ad}P)$, $\beta \in \Omega^1({\rm ad}S)$,
and $\Lambda_3 \in \Omega^3(Y^{11})$ is a closed differential form on 
$Y^{11}$.
\label{C3 invariance}
\end{proposition}

Note that in terms of differential characters, $\Lambda_3$ will be integral as
in \cite{DFM}.

\vspace{3mm}
\noindent {\bf Harmonic part of the $C$-field.}
The $C$-field has a classical harmonic part, which we now characterize.
The Bianchi identity and equation of motion for the $C$-field
in M-theory are
 \bea
 dG_4 &=&0
 \label{Bianchi}
 \\
\frac{1}{\ell_p^3} d*G_4 &=& \frac{1}{2}G_4 \wedge G_4 - I_8,
 \label{EOM}
 \eea
where $I_8$ is the one-loop term \cite{VW} \cite{DLM} 
$
I_8=\frac{1}{48}\left(p_2 - (\frac{1}{2}p_1)^2\right)
$, a polynomial
in the Pontrjagin classes $p_i$ of $Y^{11}$, $*$ 
is the Hodge duality operation in eleven dimensions,
and $\ell_p$ is the scale in the theory called the Planck constant.  
The classical (or low energy) limit given by eleven-dimensional supergravity, 
is obtained by taking $\ell_p \to 0$ and is dominated by the metric-dependent 
term. The other limit is the high energy limit probing M-theory 
and is dominated by the topological, i.e. metric-independent terms.

\vspace{3mm}
Let
$
\Delta_g^3 : \left(\Omega^3(Y^{11}), g \right)
\longrightarrow
\left(\Omega^3(Y^{11}), g \right)
$
be the Hodge Laplacian on 3-forms on the base $Y^{11}$ with respect
to the metric $g$ given by $\Delta_g^3= d~d^{{}^*}+d^{{}^*}d$, where
$d^{{}^*}$ is the adjoint operator to the de Rham differential operator 
$d$. 
Assuming $[G_4]=0$ in $H^4(Y^{11};\R)$ so that
$G_4=dC_3$, then applying the Hodge operator on (\ref{EOM})
 gives
\begin{proposition} In the Lorentz gauge, $d^{{}^*}C_3=0$, we have
\begin{enumerate}
\item $\Delta_g^3 C_3= * j_e$, where $j_e$ is the electric current associated
with the membrane given by
\(
j_e= \ell_p^3 \left(\frac{1}{2} G_4 \wedge G_4 - I_8 \right)\;.
\label{je}
\)
\item $C_3$ is harmonic if $\ell_p \to 0$ and/or there are no membranes.
\end{enumerate}
\label{delta c3}
\end{proposition}

\vspace{3mm}
The space of harmonic 3-forms on $Y^{11}$ is 
$
{\H}_g^3 (Y^{11}) : = {\rm ker} \Delta_g^3 \subset \Omega^3(Y^{11})
$.
We would like to consider harmonic 3-forms on the Spin bundle $S(Y^{11})$. 
Let
$
\Delta_{g_{\delta}}^3 : \left(\Omega^3(S(Y^{11})), g_{\delta} \right)
\to 
\left(\Omega^3(S(Y^{11})), g_{\delta} \right)
$
be the Hodge Laplacian for 3-forms on $S(Y^{11})$ with respect to the metric
$g_{\delta}$. The harmonic forms, which are in ${\rm ker} \Delta_{g_{\delta}}^3$, 
on the Spin bundle can be calculated in the adiabatic limit $\delta \to 0$. 
The expression 
for $ {\rm ker} \Delta_{0}^3:=\lim_{\delta \to 0} {\rm ker} \Delta_{g_{\delta}}^3$ was
calculated in \cite{Red}, using the spectral sequence of \cite{MaM}, further
developed in \cite{Dai} \cite{For}. Applying the results of \cite{Red} to our case  
gives
\begin{proposition}
\label{har}
\begin{enumerate}
\item When $\frac{1}{2}p_1(Y^{11}) \neq 0 \in H^4(Y^{11};\R)$, i.e. 
$[ {\rm Tr} R \wedge R] \neq 0 \in H^4(Y^{11};\R)$ then 
$
{\rm ker} \Delta_0^3 = \pi^* {\H}_g^3 (Y^{11}) \subset \Omega^3(S(Y^{11}))
$.
\item  When $\frac{1}{2}p_1(Y^{11})=0$ then 
\(
{\rm ker} \Delta_0^3 = \R \left[ CS_3(\omega) - \pi^* h \right] \oplus
 \pi^* {\H}_g^3 (Y^{11}) \subset \Omega^3(S(Y^{11}))\;,
 \label{ker delta}
\)
where $h \in \Omega^3(Y^{11})$ is the unique form such that
 $dh = {\rm Tr} R \wedge R$, $h \in d^{{}^{*}} \Omega^4(Y^{11})$\;.
\end{enumerate}
\end{proposition}
We see that for the $C$-field in M-theory we have
\begin{proposition}
\begin{enumerate} 
\item When $Y^{11}$ is a Spin manifold such that $\frac{1}{2}p_1(Y^{11}) \neq 0$ ,
 the little $c$-field 
is a harmonic form both on $Y^{11}$ and on $S(Y^{11})$. 

\item When $Y^{11}$ is a String manifold, i.e. with $\frac{1}{2}p_1(Y^{11})=0$ 
in cohomology, so that $[G_4] = a$, then 
there is a gauge in which the 3-form part of the $C$-field is that 
defining a
String class, as in the above discussion.
\end{enumerate}
\end{proposition}

\begin{remark}
The combination of forms appearing in proposition 
\ref{har} are exactly the ones also appearing in heterotic string theory. 
Indeed, the Chapline-Manton coupling is a statement about the 
String class. 
\end{remark}


\vspace{3mm}
\noindent {\bf The String class from the String condition on the target $Y^{11}$.}
From (\ref{g4 formula}), we see that when $G_4= {\rm Tr} F\wedge F$ then
$
\frac{1}{2} {\rm Tr} R \wedge R = dc
$.
At the level of cohomology this means that $\frac{1}{2}p_1(S, \omega)=0$, 
i.e. that our space admits a String structure $\S$. 
Let us form the combination $CS_3(\omega) - \pi^*c \in
\Omega^3(S)$.
Consider a choice of String structure $\S \in H^3(S(Y^{11});\Z)$. 
From 
(\ref{ker delta}), using the results in \cite{Red},
the adiabatic limit of the harmonic representative
of $\S$ is given by 
\(
[\S]_0 := {\lim}_{\delta \to 0} [\S]_{g_{\delta}}=
CS_3(\omega) - \pi^* c_3 \in
\Omega^3(S),
\)
where the form $c_3 \in \Omega^3(Y^{11})$ has the properties:
\bea
 dc_3&=& \frac{1}{2} p_1(S, \omega),
 \\
 d^{{}^*} c_3&=&0 ~~~~({\rm Lorentz~condition})\;.
\eea

\begin{observation}
Under these assumptions, the $C$-field can be identified with a String class. 
\end{observation}

\begin{remark} 
\label{cha}
Consider a change of the String structure $\S$. If the String structure is 
changed by $\xi \in H^3(Y^{11};\Z)$ then the cohomology 
class of the String structure changes, in the adiabatic limit, 
as 
\(
[\S + \pi^* \xi]_0 = [\S ]_0 + \pi^* [\xi]_g\;.
\)
Since $\xi$ is a degree three cohomology class, 
the field strength $G_4$ does not see the change in String structure.
However, at the level of the exponentiated $C$-field, i.e. at the
level of holonomy or partition function of the membrane, there will 
be an effect. See the discussion leading to observation \ref{dep st st}.
This will be formalized and considered in more detail in a future work.
\end{remark}

\subsection{The String class from the membrane.}
\label{st mem}
Consider the embedding of the membrane $\Sigma^3 \hookrightarrow Y^{11}$
with normal bundle ${\mathcal N}$. 
The fields on the membrane worldvolume include a metric $h$, 
the pullback of the $C$-field and a spinor $\psi \in \Gamma (S(Y^{11})|_{\Sigma^3})$.
Restricting $S(Y^{11})$ to $\Sigma^3$ gives the splitting
\(
S(Y^{11})|_{\Sigma^3}= S(\Sigma^3) \otimes S^{-} ({\mathcal N})
\oplus S(\Sigma^3) \otimes S^+ ({\mathcal N}),
\)
where kappa symmetry-- a spinorial gauge symmetry-- requires the 
fermions to be sections of the first 
factor \cite{BST}.
Taking the membrane as an elementary object, the exponentiated action  
will contain a factor 
$
\exp \left[ i \int_{\Sigma^3} ( C_3 + i \ell_p^{-3} {\rm vol} (h) \right]
$.
This is one part of the partition function, with another being the spinor part
given by the Pfaffian of the Dirac operator. Neither of the factors in 
\(
Z_{M2}={\rm Pfaff}(D_{S(\Sigma^3) \otimes S^{-} ({\mathcal N})}) 
\exp \left[ i \int_{\Sigma^3} ( C_3 + i \ell_p^{-3} {\rm vol} (h) \right]
\label{pf}
\)
are separately well-defined, but the product is \cite{Flux}.
Taking $\Sigma^3$ to be the boundary of a 4-manifold $B^4$ 
we get 
$
\int_{B^4} G_4
$
in place of the first factor in the exponent  in (\ref{pf}). The 
partition function is independent of the choice of 
bounding manifold $B^4$.

\vspace{3mm}
The quantization condition (\ref{shift})
for the $C$-field in M-theory was derived in \cite{Flux} by studying the 
partition function of the membrane of worldvolume $M^3$ 
embedded in spacetime $Y^{11}$. There, the manifold $M^3$ was 
assumed to be Spin. Here we notice that 
$M^3$ already admits a String structure because
$
\frac{1}{2}p_1(M^3)=0
$,
by dimension reasons. Since this is automatic, one might
wonder what is gained by assuming this extra structure. 
We will proceed with justifying this.

\vspace{3mm}
The idea is that while $M^3$ always admits a String structure, we can have more 
than one String structure. We have the following diagram
\(
\xymatrix{
K(\Z, 3) 
\ar[r]
&
B{\rm String} 
\ar[r]
\ar[d]
&
\ast
\ar[d]
\\
M^3 
\ar[u]
\ar[ur]^{\psi}
\ar[r]^{\sigma}
&
B{\rm Spin}
\ar[r]^{\lambda}
&
K(\Z, 4)\;.
}
\label{M3}
\)
Choosing a String structure $\psi$ is equivalent to trivializing 
$\lambda \circ \sigma$.  
If we fix one String structure $\psi$ then any other is classified by maps from 
$M^3$ to $K(\Z, 3)$, which is $H^3(M^3;\Z)$.  If we take 
$K(\Z, 3)={\rm B}K(\Z, 2)$ then we can say that 
the set of String structures
on $M^3$ 
is a torsor over the group of equivalence classes of gerbes on 
$M^3$. From one given (equivalence class of) String structure
we obtain for each (equivalence class of a) gerbe another
(equivalence class of a) String structure.
Notice that 
in the non-decomposible part of the $C$-field 
$
c_3 + h_3
$,
 $h_3$ is the curvature of the gerbe. It is closed as a differential form. 
We can see that there is a gerbe on the membrane worldvolume 
by taking the membrane to be of open topology and having a boundary
on the fivebrane. There is then a gerbe connection $C_3 - db_2$, where 
$b_2$ is the chiral 2-form on the M5-brane worldvolume (see e.g. \cite{AJ}).
In this case we have the exponentiated action
\(
\exp \left[ 
i \left( \int_{\partial \Sigma^3} b_2 + i \int_{\Sigma^3} \ell_p^{-3} {\rm vol}(h)
\right)
\right]\;.
\)
(See also Remark \ref{cha} and the discussion around equation (\ref{seq})).

\vspace{3mm}
Now we consider the bounding 4-dimensional space $B^4$, 
$\partial B^4=M^3$. 
Starting with a Spin $M^3$, 
there are two cases to consider, according to whether $B^4$ is 
Spin or String. 
We will make use of an approach due to David Lipsky and 
to Corbett Redden.
Let us start with the Spin case. Including $B^4$ in 
diagram (\ref{M3}) we get 
\(
\xymatrix{
K(\Z, 3) 
\ar[r]
&
B{\rm String} 
\ar[r]
\ar[d]
&
\ast
\ar[d]
\\
M^3 
\ar[u]
\ar[ur]^{\psi}
\ar[r]^{\sigma}
~ \ar@{^{(}->}[d]
&
B{\rm Spin}
\ar[r]^{\lambda}
&
K(\Z, 4)
\\
B^4
\ar[ru]^{\tau}
&&
\;.
}
\label{B4}
\)
Let $CM^3$ be the cone on $M^3$. The fact that 
\(
\xymatrix{
M^3 
\ar[r]^{\hspace{-3mm} \psi}
&
B{\rm String}
\ar[r]
\ar[d]
& 
\ast
\ar[d]
\\
& 
B{\rm Spin} 
\ar[r]
& 
K(\Z, 4)
}
\)
commutes up to homotopy means precisely that there is a strictly 
commuting diagram
\(
\xymatrix{
&
M^3 \ar[r]
\ar[d]^{i_1}
& 
\ast
\ar[ddd]
\\
M^3 
\ar[d]^{\psi}
\ar[r]^{i_0}
&
M^3 \times I
\ar[ddr]
&
\\
B{\rm String}
\ar[d]
&
&
\\
B{\rm Spin}
\ar[rr]
&&
K(\Z, 4)
}\;.
\label{new diag}
\)
Moreover, the cone is precisely the pushout
\(
\xymatrix{
M^3 
\ar[r]
\ar[d]
&
\ast
\ar[d]
\\
M^3 \times I 
\ar[r]
&
CM^3
}\;,
\)
through which, hence, diagram (\ref{new diag}) factors
\(
\xymatrix{
&
M^3 \ar[r]
\ar[d]^{i_1}
& 
\ast
\ar[d]
\\
M^3 
\ar[d]^{\psi}
\ar[r]^{i_0}
&
M^3 \times I
\ar[ddr]
\ar[r]
&
CM^3
\ar[dd]^{C\psi}
\\
B{\rm String}
\ar[d]
&
&
\\
B{\rm Spin}
\ar[rr]
&&
K(\Z, 4)
}\;.
\label{2nd new diag}
\)
Here $C\psi$ is the extension of $\psi$ 
from $M^3$ to the cone on $M^3$.
Therefore, we find that the
following diagram 
\(
\xymatrix{
M^3 
~\ar@/^2pc/[rrr]^{\lambda \circ \sigma} 
\ar[r]
\ar@{^{(}->}[d]
&
CM^3 
\ar[rr]^{C \psi}
\ar[d]^{\rho}
&&
\ast
\ar[d]
\\
B^4
\ar[r]
\ar@/_2pc/[rrr]_{\lambda \circ \tau} 
&
B^4 \bigcup_{{}_{M^3}} CM^3
\ar[rr]^{\quad \lambda(\tau, \psi)}
&&
K(\Z, 4)
}
\)
commutes.
Note that the fact that $\psi$ extends from $M^3$ to the cone of 
$M^3$, as indicated, is crucially another incarnation of the fact that
$\psi$ is homotopic to the map through the point. 
The map $\rho$ is equivalent to a String structure, and the map  
$\lambda(\tau, \psi)$ is the relative String class. Let $[B^4, M^3]$
be the relative fundamental class and $\langle ~ , ~ \rangle$ the pairing
between cohomology and homology. For this pairing we will study 
integrality and (in)dependence on the choice of $B^4$ or structures
on $B^4$.
The long exact sequence
for relative cohomology is
\bea
\cdots
\to
 H^3(M^3) 
 &\longrightarrow&
H^4(B^4, M^3)
\longrightarrow
H^4(B^4)
\longrightarrow
H^4(M^3)
\cdots
\nonumber\\
 \alpha \omega_3 
&\mapsto&
\lambda(\tau, \psi) + \partial(\alpha \omega_3)
\longmapsto
\lambda(\tau)\;,
\label{seq}
\eea
where $\omega_3$ is the volume form on $M^3$ and $\alpha$ is a 
real number. The sequence then is explained as follows.
 Given a 
choice of initial String structure on $M^3$, any other choice will 
be given by the difference with multiples of the volume form
$\omega_3$. Note that the only parameter which is varying is 
$\alpha \in \R$. 

\vspace{3mm}
Then, let $B'^4$ be another  
bounding manifold with corresponding Spin structure $\tau'$. 
Then, by the index theorem, 
\(
\int_{B^4} \frac{\lambda(\tau, \psi)}{24} -
\int_{B'^4} \frac{\lambda(\tau', \psi)}{24} = x \in \Z\;,
\)
so that
\(
\int_{B^4} \lambda(\tau, \psi) -
\int_{B'^4} \lambda (\tau', \psi) = 24 x,~~~x \in \Z\;.
\)
Now taking $e^{2\pi i}$ of both sides gives 
that the expression is an integer. So this together 
with what we explained above also gives that 
 integral does not depend on the choice of $B^4$ or structures
on $B^4$, and we have the following

\begin{proposition}
In the case when $M^3$ is String and $B^4$ is Spin, the relative pairing 
$\langle \lambda(\tau, \psi)~,~ [B^4, M^3] \rangle
= \int_{B^4} \lambda( \tau, \psi)$ is
well-defined mod 24. 
\label{spin 24}
\end{proposition}

\vspace{3mm}
Now consider the case when $B^4$ also admits a String structure,
so that $\int_{B^4} \lambda=0$. 
In this case, the question simply reduces to a statement in cobordism
of String 3-manifolds
$
\Omega_3^{\langle 8 \rangle}
 {\buildrel{\cong}\over{\longrightarrow}} \Z/24$
 defined by 
$(M^3, \psi) \longmapsto  \int_{B^4}\lambda(\tau, \psi)$.
Therefore
\begin{proposition}
In the case when both $M^3$ and $B^4$ are String manifolds, 
the relative pairing 
$\langle \lambda(\tau, \psi)~,~ [B^4, M^3] \rangle
= \int_{B^4} \lambda( \tau, \psi)$ is
well-defined mod 24. 
\label{string 24}
\end{proposition}


\vspace{3mm}
\noindent {\bf Dependence of the membrane partition function on the String structure.}
Taking (\ref{pf}) into account and the fact mentioned above
(in the proof of proposition \ref{spin 24} that 
changing the String structure of the membrane amounts to changing 
its  volume, we have for membrane worldvolumes with String 
structure

\begin{observation}
The membrane partition function depends on the choice of 
String structure on the membrane worldvolume. 
\label{dep st st}
\end{observation}

\begin{remarks}
\begin{enumerate}
\item There are nonperturbative effects, namely instantons,
 resulting from membranes wrapping 3-cycles in spacetime.
 See for example \cite{HM}.
\item In the case of string theory, the 
partition function depends crucially on the Spin structure of the string 
worldsheet. Modular invariance requires
summing over all such structures \cite{SW}. The
observations we made above then suggest that the membrane theory 
would require  a careful consideration of dependence on 
the String structure, and possibly summing over such structures. 
We hope to address this important issue elsewhere.  
\end{enumerate}
\end{remarks}

\vspace{3mm}
\noindent {\bf The framing in Chern-Simons theory.}
We can look at the dependence of the membrane partition function
on the String structure through the dependence of Chern-Simons theory on 
the choice of framing.
The partition function of  Chern-Simons theory on $M^3$ depends on
 \cite{Jones}:
$M^3$, the structure group $G$,
 the Chern-Simons coupling $k$, and a choice of framing  
$f$ of the manifold. In 
particular, the semiclassical partition
 function, while independent of the metric, it does 
 depend on the choice of framing, and different framings
 generally give different values for the partition function. 
 However, there are transformations that  map 
 the value corresponding to one  framing to the value 
 corresponding to another. 
 A framing of $M^3$ is a homotopy class of a trivialization of the
 tangent bundle $TM^3$. Given a framing $f: M^3 \to TM^3$ of $M^3$ 
 the gravitational Chern-Simons term can be defined as
 \(
 I_{M^3}(g, f)=\frac{1}{4\pi} \int_{M^3} f^* CS(\omega)\;,
 \)
 where $g$ is the metric on $M^3$, $\omega$ is the Levi-Civita
 connection on $M^3$, and the integrand is the pullback 
 via $f$ of the Chern-Simons form on $TM^3$.

\vspace{3mm}
In dimension three, there is an isomorphism 
between the String cobordism group $\Omega^{\rm String}_3$ and 
the framed cobordism group $\Omega^{\rm fr}_3$. 
Thus, the study of Chern-Simons theory with a String structure
is then equivalent to the study of Chern-Simons theory with
a framing. Therefore, it is natural to consider a String structure 
in Chern-Simons theory, and hence on the membrane worldvolume,
since the latter is essentially described by Chern-Simons theory.

\vspace{3mm}
We argue that not only is a String structure allowable, but 
is in fact {\it desirable}. This is because such a structure 
explains the framing anomaly in a very natural way. 
Under the transformation $I_{M^3} \to I_{M^3} + 2\pi s$,
where $s$ is the change in framing,  the partition function
transforms as \cite{Jones}
$Z_{M^3} \to Z_{M^3} \cdot \exp \left(2\pi is \cdot \frac{d}{24}\right)$,
for $d$ a constant related to the level of the theory. 
This factor of 24, making the partition function essentially 
a 24th root of unity, is reflection of the fact that both the String-
and the framed cobordism 
group are isomorphic in dimension three to $\Z/24$. Any Lie group $G$ 
has two canonical String structures defined by the left 
invariant framing $f_L$ and the right invariant framing $f_R$
 of the tangent bundle $TG$. For example taking $M^3$ to be
 $G=SU(2)$, there are three framings: a left framing and a right framing
 (related by orientation reversal) and a trivial framing given by 
 taking $S^3=\partial D^4$ to be the boundary of the 4-disk. 
 The invariants associated to these framings are 
 the images of $\Omega_3^{\rm String}$
 under the $\sigma$-map (the String orientation \cite{AHS})
 in $tmf^{-3}\cong \Z/24$. This map depends on the String 
 structure in an analogous way that its more classical 
cousin, the Atiyah $\alpha$-invariant, refining the $\widehat{A}$-genus 
from $\Z$ to $KO$, depends on the Spin structure. 
 From the isomorphisms
 $\pi_3 S^0 \to \pi_3 M{\rm String} \to \pi_3 tmf$ 
 \cite{Ho}
  these values are as follows
 \bea
 \Omega_3^{\rm String} &\longrightarrow & tmf^{-3}
 \nonumber\\
 \left[ SU(2), f_L \right] &\longmapsto & -\frac{1}{24}
 \nonumber\\
  \left[ SU(2), f_R \right] &\longmapsto & \frac{1}{24}
 \nonumber\\
  \left[ SU(2), \partial D^4 \right] &\longmapsto & 0\;.
 \eea
The transformation of the partition function can then be soon
more transparently using String cobordism. Thus, we get
more confirmation to observation \ref{dep st st} and, in fact, we
can also add 
\begin{observation} 
The dependence on framing of Chern-Simons theory (and hence also for 
the membrane partition function) can be seen more naturally within 
String cobordism. We thus conjecture that the membrane partition function
takes values in (twisted) $tmf$. 
 \end{observation}

Note that within Spin cobordism there would be no nontrivial expressions
in dimension three. This is because $\Omega_3^{\rm Spin}=0$ and also 
the target for the Atiyah $\alpha$-invariant, $KO_3({\rm pt})$, is also zero.
Furthermore any generalization of $\alpha$, for example to the 
Ochanine genus with target $KO_3({\rm pt})[[q]]$ would also be trivial. 

\vspace{3mm}
\noindent {\bf Further relation to generalized cohomology.}
There is further a connection of the membrane partition function 
to generalized cohomology
from another angle as follows. 
The Stiefel-Whitney class 
$w_4$ is the mod 2 reduction of $\lambda=\frac{1}{2}p_1$. This implies that 
$\lambda$ is even
if and only if $w_4=0$. The latter is in fact an orientation condition in real Morava
$E$-theory $EO(2)$ \cite{KS1} (worked out there for a different but related purpose). 
Therefore, in order to remove the ambiguity in the
quantization, we require $EO(2)$-orientation.

 \begin{proposition}
 The membrane partition function is well-defined in 
 $EO(2)$-theory. In particular, $G_4$ is an integral class
 when the underling spacetime is $EO(2)$-oriented. 
 \end{proposition}
 
 \begin{remarks}
 \begin{enumerate}
\item This means that eleven-dimensional 
 spacetime backgrounds in M-theory with no fluxes 
 should be $EO(2)$-oriented.
\item Note that $EO(2)$ is closely related to the theory 
$EO_2$  of Hopkins and Miller, which in turn is closely 
related to $TMF$.
\item The M-theoretic partition function via $E_8$ gauge theory
of \cite{DMW} is 
considered for the String case in \cite{MO8}.
\end{enumerate}
 \end{remarks}

\vspace{3mm}
\noindent {\bf The case when the Spin bundle is trivial.}
Consider the case when $S(Y^{11})$ is trivial as a principal
bundle. This then means that there is a global section. Such 
a section $s: Y^{11} \to S(Y^{11})$ gives an isomorphism 
of principal bundles
\(
S(Y^{11}) \buildrel  s \over{ \cong}
Y^{11} \times {\rm Spin}(11)\;.
\label{section}
\)
Using the fact that $H^i\left( {\rm Spin}(11);\Z \right)=0$ for 
$i=1,2$, the K\"unneth formula gives an induced isomorphism 
on integral cohomology 
\bea
H^3(S;\Z) &\cong & H^3(Y^{11};\Z) \oplus H^3({\rm Spin}(11);\Z)
\nonumber\\
\S &\longleftrightarrow & (0~,~ 1_{\rm Spin})\;.
\eea
The String structure is then determined by $s$ and is an element 
$\S \in H^3(S;\Z)$ which corresponds to the pullback of the
generator $1_{\rm Spin} \in H^3({\rm Spin}(11);\Z)$.

\vspace{3mm}
Note that the differential $d$ acting on $p^*CS_3(\omega)$
is $\frac{1}{2}p_1(\omega)$, which is the same as
$d c$. This means that 
$[p^* CS_3(\omega) - c] \in H^3(Y^{11};\R)$. 
Consider the membrane worldvolume $M^3$, taken as a 3-cycle
$X \in {\rm Map}(M^3, Y^{11})$. Then 
$
p(X) \subset Y^{11} \times \{ {\rm pt}\} \subset S
$ 
under the isomorphism (\ref{section}) 
induced by the global section $s$. The triviality of the bundle
implies that any class $\xi \in H^3(S;\Z)$ is zero when evaluated on 
3-cycles $\Sigma^3$ in $Y^{11}\subset S$,
$
\langle \xi , [s(\Sigma^3)] \rangle,
$, 
where $[s(\Sigma^3)]$ is the fundamental class of the 3-cycle $s(\Sigma^3)$ 
in $S(Y^{11})$. Then, in real cohomology
\(
0 = \int_{s(\Sigma^3)} [\xi]_0 
= \int_{s(\Sigma^3)} \left( CS_3(\omega) - \pi^*c \right)
= \int_{\Sigma^3} \left( s^* CS_3(\omega) - c \right)\;.
\)
This holds for an arbitrary 3-cycle $\Sigma^3$ so that 
$
[s^*CS_3(\omega) - c] =0 \in H^3(Y^{11};\R)
$.
Now if $d^{{}^*}(s^*CS_3(\omega))=0$, the since 
$d^{{}^*} c=0$, then 
$s^*CS_3(\omega) - c$ is harmonic, and
hence zero, so
$
s^* CS_3(\omega) = c
$.
Therefore, in this case, following the general construction 
\cite{Red} \cite{Red2}, $s^* CS_3(\omega)$ and $c$ are 
equal to the coexact and harmonic components, respectively.

\vspace{3mm}
\noindent {\bf The $C$-field as a Chern-Simons 2-gerbe.}
In \cite{CJMSW} a Chern-Simons bundle 2-gerbe is constructed, 
realizing differential geometrically the Cheeger-Simons
invariant \cite{CS}.  This is done by introducing a lifting
to the level of bundle gerbes of the transgression map from 
$H^4(BG; \Z)$ to $H^3(G; \Z)$. A similar construction 
is given in \cite{AJ}. Both groups of interest, ${\rm Spin}(n)$
and $E_8$ are simply-connected, a fact that removes some 
subtleties from the discussion.

\vspace{3mm}
For any integral cohomology class in $H^3({\rm Spin}(n); \Z)$, there is a 
unique stable equivalence
class of bundle gerbes (\cite{Mu} \cite{MuS}) whose Dixmier-Douady class is the given
degree three integral cohomology class. 
Geometrically $H^4(B{\rm Spin}(n); \Z)$ 
can be regarded as stable equivalence classes of bundle 2-gerbes 
over $B{\rm Spin}(n)$, whose induced bundle gerbe
over ${\rm Spin}(n)$ has a certain multiplicative structure.
More precisely, given a bundle gerbe $\mathcal{G}$ over ${\rm Spin}(n)$, 
${\mathcal G}$ is multiplicative if and only if its
Dixmier-Douady class is transgressive, i.e., in the image of the transgression map
$\tau : H^4(B {\rm Spin}(n); \Z) \to H^3 ({\rm Spin}(n);\Z)$ \cite{CJMSW}.



\vspace{3mm}
Consider a principal ${\rm Spin}(n)$-bundle $P$ 
with connection $A$ on a manifold $M$. 
For the Chern-Simons gauge theory canonically defined by a class
in
$H^4(B {\rm Spin}(n); \Z)$, 
there is a Chern-Simons 
bundle 2-gerbe $Q(P,A)$ associated with the
$P$ which is  
defined to be the pullback of the universal Chern-Simons bundle 
2-gerbe by the classifying map $f$ of $(P,A)$ \cite{CJMSW}.

\vspace{3mm}
With the canonical isomorphism between the Deligne cohomology
and Cheeger-Simons cohomology, the Chern-Simons bundle 2-gerbe $Q(P,A)$ 
is equivalent \cite{CJMSW} in Deligne cohomology to the Cheeger-Simons 
invariant $S(P,A) \in \check{H}^3(M, U(1))$, which is 
the differential character that can be 
associated with each principal $G$-bundle $P$ with 
connection $A$ \cite{CS}. A bundle gerbe version of the discussion of 
the invariance of the $C$-field in section \ref{c as string} is provided by 
the following

\begin{proposition}
[\cite{AJ} (also \cite{CJMSW})]
A Chern-Simons 2-gerbe is contained in the data of the $C$-field. 
\label{AJ prop}
\end{proposition}

\vspace{3mm}
\noindent{\bf The metric torsion part of $C$-field.}
Connections with torsion come in various classes. Especially interesting is the
case when the torsion is totally antisymmetric. In this case the new connection is
metric and geodesic-preserving, and the Killing vector fields coincide with the
Riemannian Killing vector fields. 
Connections with torsion arise in  
in eleven-dimensional supergravity \cite{Englert} and in heterotic string theory
and type I \cite{Str}. In the first case an ansatz is taken such that the little
$c$-field is proportional to the torsion tensor $T \in \Omega^3(Y^{11})$,
most prominently when $Y^{11}$ is an $S^7$ bundle over 4-dimensional 
anti-de Sitter space AdS$_4$, in which case the torsion is parallelizing 
-- see \cite{DNP}. In the second case, the $H$-field in heterotic string theory
acts as torsion, which is important for compactification to lower dimensions
\cite{Str}.


\section{The M-theory Dual $C$-Field and Fivebrane Structures}

\subsection{The $E_8$ model for the dual of the $C$-field}

In \cite{DFM}  the electric charge associated with the $C$ field is studied.
From the nonlinear equation of motion (\ref{EOM})
of the $C$-field, the induced electric charge that is given by 
the cohomology class 
\(
[\frac{1}{2}G^2-I_8]_{DR}\in H^8(Y, \R)\,. \label{drcharge}
\)
In \cite{DFM} the integral lift of (\ref{drcharge}) 
is studied and denoted $\Theta_Y(\check C)$ (and also $\Theta_Y(a)$),
where $\check C=(A,c)$.

\vspace{3mm}
A tubular neighbourhood of the M5-brane worldvolume $V$ in Y
is diffeomorphic to 
the total space of the normal bundle $N\rightarrow V$.
 Let $X=S_r(N)$ be the 10-dimensional sphere bundle of
radius $r$, so that the fibers of $X {\stackrel{\pi}{\to}}  V$ are 4-spheres. 
An 11-dimensional manifold $Y_r$ with boundary $X$ is
then constructed by removing from $Y$ the disc bundle $D_r(N)$ of radius $r$.
$Y_r = Y-D_r(N)$, the bulk manifold, is the complement of the
tubular neighbourhood $D_r(N)$). There are two path integrals or wavefunctions:
\begin{enumerate}
\item The bulk $C$ field path integral 
$\Psi_{{\!bulk}}(\check C_X)\sim \int
\exp[G\wedge* G] \Phi(\check C_{Y_r})$ where the integral is over all
equivalence classes of $\check C_{Y_r}$ fields that on the boundary assume
the fixed value $\check C_X$. This wavefunction is a
section of a line bundle $\mathcal L$ on the space of $\check C_X$ fields. 

\item The M5-brane partition
function $\Psi_{M5}(\check C_V)$, which depends on the $\check C$ field 
on an infinitesimally small ($r\rightarrow 0$) tubular 
neighbourhood of the M5-brane.
\end{enumerate}

\vspace{3mm}
In general, $\Psi_{{\!bulk}}\Psi_{M5}$ is not gauge invariant
and therefore it is a section of a line bundle. However, one can consider 
a new partition function 
$\Psi'_{M5}$ that is obtained from multiple M5-branes stacked on top of one another
instead of just
a single M5-brane. This stack gives rise to a twisted gerbe 
on $V$ as follows \cite{AJ}. 
In order for
$\Psi_{{\!bulk}}\Psi_{M5}'$ to be well defined, 
the twisted gerbe has to satisfy 
\(\label{gentwist}
[CS({\pi_*(\Theta_X)})]-[\vartheta_{ijkl},0,0,0]=[
{\mathbf D}_H]+[1,0,0,C_V]\,,
\)
where, 
$CS({\pi_*(\Theta_X)})$ is the Chern-Simons 2-gerbe associated
with $\pi_*(\Theta_X)$ and a choice 
of connection on the $E_8$ bundle with first Pontryagin class 
${\pi_*(\Theta_X)}$ (all other 2-gerbes 
differ by a global 3-form),
while $[\vartheta_{ijkl},0,0,0]$ is the 2-gerbe class
associated with the torsion class $\theta$ on $V$,
 $\beta(\vartheta)=\theta$, and $[1,0,0,C_V]$
is the trivial Deligne class associated with the global 3-form $C_V$.
In particular (\ref{gentwist}) implies
\(
\pi_*(\Theta_X)-\theta=
\xi_{{\mathbf D}_H}\,,
\)
where the RHS is the characteristic class of the lifting 2-gerbe.

\subsection{Twisted Fivebrane structures}

Fivebrane structures are obtained by lifting String structures as follows.
\begin{definition}[{\cite{SSS1}}{\cite{SSS2}}]
 An $n$-dimensional manifold $X$ has a
 \emph{ Fivebrane structure} 
 if the classifying
 map $X \to BO(n)$ of the tangent bundle $T X$ lifts to the
 classifying space
 $B{\rm Fivebrane} :=B\mathrm{O}\langle 9 \rangle$.
\(
  \raisebox{20pt}{
  \xymatrix{
    &
     \mathrm{B O} (n) \langle 9 \rangle 
     \ar[d]
     \\
     M
     \ar[r]_f
     \ar@{..>}[ur]^{\hat f}
     &
     \mathrm{B O}(n)~.
  }
  }
\)
\end{definition}

\begin{theorem}[{\cite{SSS2}}]
\begin{enumerate}
\item The obstruction to lifting a String structure on $X$
    to a Fivebrane structure on $X$ is the fractional
    second Pontrjagin class
    $ \frac{1}{6}p_2(T X)$.\\
    \item  The set of lifts, i.e. the set of Fivebrane structures for a fixed String structure
in the real case,  or the set of $BU\langle 9\rangle$ structures for a fixed 
$BU\langle 7 \rangle$ structure in the complex case, is a 
torsor for a quotient  
of the seventh integral cohomology group $H^7(X; \Z)$. 
\end{enumerate}
\end{theorem}

\vspace{3mm}
\noindent{\bf Twisted Fivebrane structures.}
 Twisted cohomology refers to a cohomology which is
defined  in terms of a  twisted differential on ordinary differenial forms.
The twisted de Rham complex $\Omega^{\bullet} (X, d_{H_{2i+1}})$
means the usual de Rham complex but with the differential $d$
replaced by $d_{H_{2i+1}}:= d + H_{2i+1}\wedge$,
 which squares to zero by virtue of the
Bianchi identity for $H_{2i+1}$, provided that $H_{2i+1}$ is closed.
The form on which this twisted differential acts would 
involve components that are $2i$ form-degrees apart, e.g. of the form
$F=\sum_{n=0}^m F_{k+2in}$. The main case considered in \cite{7twist},
and which corresponds to heterotic string theory,
corresponds to $i=3$, $k=2$, and $m=1$.
There it 
was observed that, at the rational level, the
(abelianized) field equation and Bianchi identity in heterotic string
theory can be combined into an equation given by a degree seven
twisted cohomology. This is discussed further in section
\ref{duality sym twist}. In \cite{7twist} it
was also 
proposed that such a differential should correspond to a twist of 
what was there called a higher String structure and later in 
\cite{SSS2} was defined, studied in detail and given the name 
Fivebrane structure. 
\begin{definition}[ \cite{SSS3}]
  A $\beta$-twisted Fivebrane structure
  on a brane $\iota : M \to X$
  with String structure classifying map $f : M \to BO(n) \langle 8 \rangle$ is a
  cocycle $\beta : X \to K(\mathbb{Z}, 8)$ and a map $c: BO(n) \langle  8 \rangle 
  \to K(\Z, 8)$ such that there is a homotopy $\eta$
  between  $BO(n) \langle  8 \rangle$ and $X$, 
  $$  
  \raisebox{20pt}{
    \xymatrix{
       M
       \ar[rr]^f_>{\ }="s"
       \ar[d]_\iota
       &&
       B O (n)\langle  8 \rangle 
       \ar[d]^{c}
       \\
       X
      \ar[rr]_\beta^<{\ }="t"
       &&
       K(\mathbb{Z},8)
            \ar@{=>}^\eta "s"; "t"
   }
    }\;.
    $$
   \label{2tX}
\end{definition}
The fractional class
 $\frac{1}{48}p_2$ show up in physics. 
As an application, which was originally the motivation:
\begin{theorem} [{\cite{SSS3}}]
\begin{enumerate}
\item The dual formula for the Green-Schwarz anomaly cancellation condition
on a String 10-manifold $M$ is the obstruction to defining
a twisted Fivebrane structure,
with the twist given by $ch_4(E)$, where $E$ is the gauge 
bundle with structure group $E_8 \times E_8$.
\item The  integral class in M-theory dual to $G_4$ defines an obstruction to
twisted Fivebrane structure, which
 is the obstruction to having a well-defined partition function for the
 M-fivebrane.
\end{enumerate}
\end{theorem}
\vspace{3mm}
The precise description actually requires some refinement of the Fivebrane structure
to account for appearance of fractional classes such as $\frac{1}{48}p_2$ rather than 
$\frac{1}{6}p_2$. This is called  
$\mathcal{F}^{\langle 9 \rangle}$-structure in \cite{SSS3}.
Also, a corresponding ${\mathfrak{fivebrane }}$ Lie 6-algebra
is defined in \cite{SSS3}, where the description of the twist in 
terms of $L_{\infty}$-algebras is also given. 

\vspace{3mm}
\noindent {\bf Higher bundles:} The following applications are of interest
(see \cite{SSS1}).\\
\noindent {\bf 1.}  Chern-Simons 3-forms arise as local connection data
on 3-bundles with connection which arise as the obstruction to lifts of
ordinary bundles to the corresponding String 2-bundles and 
are shown to be 
governed by the String Lie 2-algebra. For $\gg$ an ordinary semisimple 
Lie algebra and $\mu$ its 
canonical 3-cocycle, 
 the obstruction to lifting a $\gg$-bundle to a String 2-bundle
is a Chern-Simons 3-bundle with characteristic class the Pontrjagin class
 of the original bundle.  \\
\noindent {\bf 2.}  The formalism immediately allows the generalization of this situation to 
higher degrees. Indeed, certain 7-dimensional generalizations of 
Chern-Simons 3-bundles obstruct the lift of ordinary bundles to certain 6-bundles
governed by the Fivebrane Lie 6-algebra. The latter correspond what was defined in
\cite{SSS1} as the {\it Fivebrane structure}, for which the 
degree seven NS field $H_7$ plays the role
that the degree three dual NS field $H_3$ plays for the $n=2$ case. 
  Using the
  7-cocycle on $\mathfrak{so}(n)$,  lifts through 
  extensions by a Lie 6-algebra, defined as the Fivebrane Lie 
  6-algebra, is obtained. Accordingly, Fivebrane structures 
  on String structures are indeed obstructed by
  the second Pontrjagin class.  
\subsection{Harmonic part of the dual of the $C$-field}
Consider the equation of motion (\ref{EOM}). Set $C_7:=*G_4$, so that 
we get
\(
dC_7= j_e,
\label{dc7}
\)
where $j_e$ is the electric current given in (\ref{je}). 
Let 
\(
 \Delta_g^7 : \left(\Omega^7(Y^{11}), g \right)
\longrightarrow
\left(\Omega^7(Y^{11}), g \right)
\)
be the Hodge Laplacian on 7-forms on $Y^{11}$ with respect to the metric 
$g$. Taking $d^{{}^*}$ of
both sides of equation (\ref{dc7}) we get 
\begin{proposition}  In the Lorentz gauge, $d^{{}^*}C_7=0$, we have
\begin{enumerate}
\item $\Delta_g^7 C_7= * j_m$, where $j_m$ is the fivebrane magnetic current, related to the
electric current $j_e$ (\ref{je}) by
\(
j_m=d(*j_e)\;. 
\label{jm}
\)
\item $C_7$ is harmonic if $\ell_p \to 0$ and/or there are no fivebranes.
\end{enumerate}
\label{delta c7}
\end{proposition}
\noindent This is the degree seven analog of proposition \ref{delta c3}.
The geometry of the de Rham representatives of the 
Fivebrane class for the dual $C$-field is considered in \cite{RS}. 
\subsection{The integral lift}
\label{int lift}
The quantization law (\ref{shift}) on the $C$-field leads to an 
integral lift of the electric charge, the left hand side of the equation of
motion (\ref{EOM}), as mandated by Dirac quantization. The lift is \cite{DFM}
\bea
[G_8]&=& \frac{1}{2}\left( a - \frac{1}{2}\lambda \right) 
\left( a - \frac{1}{2} \lambda \right)
+ I_8
\nonumber\\
&=& \frac{1}{2} a( a - \lambda) + 30 \widehat{A}_8\;.
\label{G8}
\eea
\begin{proposition} Properties of $[G_8]$.
\begin{enumerate}
\item $[G_8]$ and $[G_4]$ obey the multiplicative structure on $K$Spin.
\item $G_8$ defines an obstruction to having a (twisted) 
Fivebrane structure.
\end{enumerate}
\end{proposition}
The first part is proved in \cite{KSpin} and the second part in \cite{SSS3}.
Let us expand a bit on the first point. The quadratic refinement
defined in \cite{DFM} is encoded in the multiplicative structure in
the K-theory for Spin bundles.
Starting with a real unoriented
bundle $\xi$, the condition $w_1(\xi)=0$ turns $\xi$ into an
oriented bundle, and the condition $w_2(\xi)=0$ further makes $\xi$
a Spin bundle. Obviously then, a real $O$-bundle becomes a Spin
bundle when $W=0$, and so the kernel of $W$ is the reduced group
 $\widetilde{K{\rm Spin}}(X)$. Thus $W$
fits into the exact sequence \cite{Li}
\( 0 \longrightarrow
{\widetilde{K{\rm Spin}}}(X)=\ker W \longrightarrow
{\widetilde{KO}}(X) {\buildrel{W} \over {\longrightarrow}}
H^1(X;\Z_2) \times H^2(X;\Z_2). \)

\vspace{3mm}
Among the properties of this class proved in \cite{DFM} is that it is
a quadratic refinement of the cup product of two degree four classes
$a_1$ and $a_2$
\(
\Theta(a_1+a_2) + \Theta(0)= \Theta(a_1)  + \Theta(a_2)+ a_1 \cup a_2.
\label{ref}
\)
We would like to look at this from the point of view of the structure on
the product of the cohomology groups $H^4(~\cdot~;\Z) \times H^8(~\cdot~;\Z)$.
For this we consider the two classes $a$ and $\Theta(a)$ as a pair
$(a,\Theta(a))$ in $H^4(~\cdot~;\Z) \times H^8(~\cdot~;\Z)$. Then the linearity
of the addition of the degree four classes $a$ and the quadratic
refinement property (\ref{ref}) of $\Theta(a)$ can both be written in
one expression in the product $H^4(~\cdot~;\Z) \times H^8(~\cdot~;\Z)$, which
makes use of the ring structure, namely
\(
\left( a_1,\Theta(a_1) \right) + \left( a_2,\Theta(a_2) \right)=
\left( a_1 + a_2, \Theta(a_1) + \Theta(a_2) + a_1 \cup a_2 \right).
\label{t2}
\)
The second entry on the RHS is just $\Theta(a_1+a_2) - \Theta(0)$,
and so it encodes the property (\ref{ref}).

\vspace{3mm}
We can define the shifted class $\Theta^0(a)$ as the difference
$\Theta(a) - \Theta(0)$, so that (\ref{t2}) is
replaced by
\(
\left( a_1,\Theta^0(a_1) \right) + \left( a_2,\Theta^0(a_2) \right) =
\left( a_1 + a_2, \Theta^0(a_1+a_2)\right),
\label{t3}
\)
corresponding to the special case
\(
\Theta^0(a_1+a_2) = \Theta^0(a_1)  + \Theta^0(a_2)+ a_1 \cup a_2.
\)
This is then just a realization of the multiplication law on
$H^4(~\cdot~;\Z) \times H^8(~\cdot~;\Z)$ which, for $(a,b)$ in the product group,
is 
\(
(a_1, b_1) + (a_2, b_2)=(a_1 + a_2, b_1 + b_2 + a_1 \cup a_2).
\label{quad}
\)
Note that in order to get this law we had to use the modified
eight-class $\Theta^0(a)$, or alternatively discard $\Theta (0)=30
{\widehat{A}}_8$.

\vspace{3mm}
We now make the connection to Spin K-theory.
Similarly to the case of other kinds of bundles, e.g. complex or real,
one can get a Grothendieck group of isomorphism classes of Spin
bundles up to
equivalence. The reduced $K{\rm Spin}$ group of a
topological space can be defined as
$\widetilde{K{\rm Spin}}(X)=[X,B{\rm Spin}]$.
For the case of $B{\rm Spin}$, we will be interested in relating Spin
K-theory to cohomology of degrees 4 and 8. Such a homomorphism of
abelian groups
\( Q_X:{\widetilde{K{\rm Spin}}}(X) \rightarrow H^4(X;\Z) \times
H^8(X;\Z)
\label{ch}
\)
is defined by \cite{Li}
$Q_X\left( Q_1(\xi),Q_2(\xi)\right)$ for $\xi \in
{\widetilde{K{\rm Spin}}}(X)$, where $Q_1$ and $Q_2$ are the Spin
characteristic classes of \cite{Th}.
For two bundles $\xi$ and $\gamma$ in
${\widetilde{K{\rm Spin}}}(X)$, and for $k\leq 3$,
\( Q_k(\xi \oplus \gamma)=\sum_{i+j=k}Q_i(\xi) \cup Q_j(\gamma).
\label{add}
\)

\begin{remark}
The above multiplicative structure is a $\Z$-analog (or $4k$-analog)
 of the
$\Z_2$-structure in the case of $KO$-theory. 
Given a topological space $X$, let ${\widetilde{KO}}(X)$
be the reduced $KO$ group for $X$ and let \( W~:~
{\widetilde{KO}}(X)\longrightarrow H^1(X;\Z_2) \times H^2(X;\Z_2)
\label{KO} \) be the map $W(\xi)=(w_1(\xi),w_2(\xi))$, where
$w_i(\xi)$ denotes the $i$-th Stiefel-Whitney class of $\xi \in
{\widetilde{KO}}(X)$. There is a group structure on $H^1(X;\Z_2)
\times H^2(X;\Z_2)$ making $W$ a homomorphism, i.e. a map that
preserves the group structure.
\end{remark}

\subsection{Invariance of the dual $C$-field}
From (\ref{G8}) we see that we can write the dual $C$-field 
at the level of differential forms as
\(
G_8=\frac{1}{2}G_4 \wedge G_4 - I_8 + dc_7\;.
\label{G8 rational}
\)
\subsubsection{Case I: Trivial cohomology, no one-loop term}
In this case we have $\frac{1}{2}p_1=0=\frac{1}{6}p_2$. 
From (\ref{G8 rational}) we have 
$
dC_7= \frac{1}{2} G_4 \wedge G_4 + dc_7$.
When $C_3=CS_3(A)$ we get
$
C_7= \frac{1}{2} CS_3(A) \wedge G_4$.

\subsubsection{Case II: $\frac{1}{2}p_1=0$} 
\vspace{3mm}
The invariance of $G_8$ will include the invariance of the three terms
in (\ref{G8 rational}), hence invariance of $G_4$, of the
Pontrjagin classes of $S(Y^{11})$, and of the differential form
$C_7 \in \Omega^7(Y^{11})$. Therefore, we have
\begin{proposition}
The dual field $G_8$ is invariant under the following transformations:
\begin{enumerate}
\item The invariances of the $C$-field from proposition \ref{C3 invariance},
\item $C_7 \mapsto C_7 + \Lambda_7$,
\item $CS_7(\omega) \mapsto CS_7(\omega) + \lambda_7$,
\end{enumerate}
where $\Lambda_7$ and $\lambda_7$ are closed differential forms in 
$\Omega^7(Y^{11})$.
\label{C7 invariance}
\end{proposition}
As in the case for the $C$-field, the expressions above when recast in terms 
of differential characters will result in requiring $\Lambda_7$ and $\lambda_7$
to be closed integral forms, i.e. to be in $\Omega^7_{\Z}(Y^{11})$.

\vspace{3mm}
\noindent{\bf The Fivebrane class.}
\begin{remark}
The Fivebrane class can be discussed in a manner that is very similar 
to that of the String class. 
The change of Fivebrane structure ${\mathcal F}$ can be 
seen from the M-theory fivebrane, similarly to the way the 
String structure ${\mathcal S}$ is seen from the membrane.
\end{remark}

Here we assume that spacetime is ten-dimensional. We 
relate the dual $H$-field $H_7$ to the degree seven Chern-Simons
form $CS_7$ above via a class $C_7$ that we introduce.  
Consider the principal bundle 
\(
{\rm String}(n) \to {\rm STRING}(M) \to M,
\label{sty}
\)
where ${\rm STRING}(M)$ is the total space of the ${\rm String}(n)$ 
bundle on our ten-dimensional spacetime $M$ corresponding to heterotic
string theory. Recall that this is one of the two bundles in that theory, namely
the one obtained from the lift of the tangent bundle (non-gauge one).
We can build a class $C_7$ out of the Chern-Simons 7-form 
$CS_7(A)$ on ${\rm STRING}(X)$ for the given connection 1-form $A$ 
on the gauge
bundle as
\(
C_7 = CS_7(A) - \pi^* H_7,
\)
with $dCS_7= \pi^*\left( p_2(M) \right)$ and such that
$\int_{\Sigma_7} C_7=1$ where $\Sigma_7$ is a fundamental 
7-cycle in the fiber. Similarly to the Spin case \cite{ASi}, we have

\begin{proposition} If $X$ is 6-connected,
$C_7$ represents a generator of $H^7({\rm STRING}(M); \Z)$.
\end{proposition}
\proof
$C_7$ is closed since
\bea
dC_7 &=& dCS_7(A) - d\pi^*H_7
\nonumber\\
&=&\pi^*\left( p_2(M)\right) - \pi^* dH_7
\nonumber\\
&=& \pi^* \left(  p_2(M) - dH_7 \right) 
\nonumber\\
&=&0.
\eea
Consider the homotopy exact sequence 
\(
\cdots \to \pi_i({\rm String}(n)) \to \pi_i({\rm STRING}(M)) \to \pi_i(M) 
\to \pi_{i-1}({\rm String}(n)) \to \cdots
\)
corresponding to the bundle (\ref{sty}). 
Assuming $\pi_i(M)=0$ for $i \leq 6$,  
the sequence gives that
$
\pi_7({\rm STRING}(M)) \cong \pi_7({\rm String}(n))=\Z
$.
 Therefore, $C_7$ represents a generator 
of $H^7\left( {\rm STRING}(M) ,Z \right)$. 
\endofproof

\begin{remark} If $M$ is ten-dimensional and is $6$-connected 
then it is topologically the ten-dimensional sphere. This follows from Poincar\'e 
duality and the Poincar\'e conjecture in ten dimensions. The first 
gives that the cohomology groups in degrees 8 and 9 are the same
as those in degrees 2 and 1, respectively, and hence are zero 
by $6$-connectedness. $M$, having cohomology groups $\Z$
in degrees $0$ and $10$, is a homology ten-sphere. However, 
by the Poincar\'e conjecture, which is a theorem in ten dimensions,
$M$ must be the sphere $S^{10}$ itself. 
\end{remark}

\subsection{Higher differential characters}
We consider the
space ${\rm Map}(Z, M)$, with $Z$ of dimension six. 
The space ${\rm Map}(M, {\rm STRING}(M))$ is a bundle over 
${\rm Map}(Z, M)$ with structure `group' ${\rm Map}(M, {\rm String}(10))$.
Let 
\(
\widetilde{ev} : Z \times {\rm Map}(Z, {\rm STRING}(M)) \to {\rm STRING}(M)
\)
be the evaluation map.

\begin{proposition}
There exists
a Cheeger-Simons differential 6-character $B_6$ with
$dB_6=C_7$ and such that
$\widetilde{ev}^* B_6$ exponentiates to a differential 0-character on 
${\rm Map}(Z, {\rm STRING}(M))$ with values in $U(1)$.
\end{proposition}

\proof 
This is analogous to \cite{ASi} where the degree two case
is established. Our case corresponds to replacing $B_2$ with $B_6$
and the spin condition with the String condition. For a trivial map
$\widetilde{\Phi}_0 : Z \to P_0 \in {\rm STRING}(M)$, let $\gamma$ be 
a path from $\widetilde{\Phi}_0$
to $\widetilde{\Phi} \in {\rm Map}(Z, {\rm STRING}(M))$ so that $\gamma$
maps the interval $[0,1]$ times $Z$ to ${\rm STRING}(M)$ with 
$\gamma(1)=\widetilde{\Phi}$. This path exists since 
$\pi_i({\rm STRING}(M))=0$ for $i\leq 6.$

The function $\widetilde{ev}^*B_6$ will be 
given by
\(
\exp \left( 2\pi i \int_{\gamma([0,1] \times Z)}C_7 \right)
=\exp \left(  2 \pi i \int_{\gamma([0,1])} \alpha_1 \right)
\label{evb6}
\)
with $\alpha_1$ is the one-form $\int_{Z} \widetilde{ev}^* C_7$. 
The function is independent of the path
$\gamma$: if $\gamma_1$ is another path then $\gamma_1^{-1}\gamma$
is a map of $S^1 \times Z \to P_0$ and
\(
\int_{(\gamma_1^{-1}\gamma)(S^1 \times Z)} C_7
=\int_{\gamma(S^1 \times Z)} C_7 - 
\int_{\gamma_1(S^1 \times Z)} C_7 \in \Z,
\)
i.e. $\alpha_1$ represents an integral 1-cocycle. 
Since $C_7$ represents an element of $H^7({\rm STRING}(M);\Z)$, there
exists a Cheeger-Simons differential 6-character $B_6$ with 
$dB_6=C_7$. Then $\int_{Z} \widetilde{ev}^*B_6$ exponentiates to
a differential  
$0$-character on ${\rm Map}(Z, {\rm STRING}(M))$ with values in $S^1$.
\endofproof

\begin{proposition}
The curvature of the $0$-character $\widetilde{ev}^*B_6$
is the one-form $\alpha_1$, which can be
interpreted as a flat connection on a circle bundle over ${\rm Map}(Z, {\rm STRING}(M))$. 
\end{proposition}

\begin{remarks}
{\bf 1.} We consider the effect of gauge transformation on $H_7$. 
Changing $H_7$ to $H_7+ d \beta_6$ leads to a shift in $B_6$ as $B_6 + \beta_6$.

{\bf 2.} The topological term 
\(
I= \int_{Z} dx^1 \wedge \cdots \wedge dx^6 B_{\mu_1 \cdots \mu_6}
\partial_1 X^{\mu_1} \cdots \partial_6 X^{\mu_6}
\)
can be interpreted as $\log \int_{Z} \widetilde{ev}^* B_6$, in analogy with
the case of the string.
 
 {\bf 3.} Similar results can be extrapolated to eleven dimensions.
\end{remarks}

    \subsection{Mapping  Space Description}
    \label{Configuration spaces}
   
  As we recalled in the introduction, String structures on a space $M$ is 
  related to Spin structure on the loop space $LM$. Is there a corresponding 
statement  about Fivebrane structures? We do not fully answer this 
question here but we do give possible scenarios.        
     
     \vspace{3mm}     
        In the case of String structure, the string class on the loop space $LM$
    is obtained by pull-back of the second Chern character, $ch_2(E)$, of 
    a bundle $E$ on $M$ to give a bundle $LE$ on $LM$ via the evaluation map 
   \(
   ev : S^1 \times LM \to M,
   \label{ev1}
   \)
   so that the String class is
   $
   \int_{S^1} ev^* ch_2(E)
   $,
   where $\int_{S^1} : H^*(S^1 \times LX; \C) \to H^{*-1}(LM; \C)$ 
   is the integration along $S^1$.
   
 
 \vspace{3mm}
   Now the idea is to generalize this to the case of the Fivebrane structure. 
   We see two directions for doing so:
   \begin{enumerate}
   \item Replace the second Chern character $ch_2$ by a higher degree Chern character
  $ch_{q}$, $q>2$, and keeping the same evaluation map (\ref{ev1}). This will yield
  higher degree analogs of the String class, but {\it still} on the loop space $LM$.
  
  \item In addition, replace the circle $S^1$ by a higher-dimensional space $Y$ so 
  that the loop space $LM=[S^1, M]$ is replaced by a higher degree generalization
 ${\rm Map}[Y,M]$, the space of maps from $Y$ to $M$.  
 \end{enumerate}

\vspace{3mm}
 We start with the first. Here, for a vector bundle $E$ over $M$,
 one has the higher degree analogs of the String 
 class as 
  \(
 C^q(LE)=- (2 \pi i)^{q+1} q! \int_{S^1} ev^* ch_{q+1}(E),
  \label{stp}
  \)
  which gives a class of degree $2(q+1)-1=2q-1$ on $LM$. Such a generalization of the
  usual String structure has been defined in \cite{As}. 
  Kuribayashi in \cite{Ku} finds fairly special conditions under which it is still true
 that $C^p(LE)=0$ if and only if $ch_{p+1}(E)=0$, namely when $H^*(M; \R)$ is 
 a tensor product of truncated polynomial and exterior algebras. This generalizes 
 McLaughlin's result \cite{Mc} in the case $p=1$ for the usual String structure, 
 where $\frac{1}{2}p_1(P)=0$ implies
a String structure on a bundle $P$ on $M$ only when $\pi_2(M)=0$.

  \vspace{3mm}
  An example of the higher classes (\ref{stp}) 
  is the first term in the anomaly
  polynomial 
  \(
  dH_7=2\pi \left[  {\rm ch}_2(A) - \frac{1}{48} p_1(\omega) {\rm ch}_2(A)
  + \frac{1}{64} p_1(\omega)^2 - \frac{1}{48} p_2(\omega) \right],
  \)
  where $A$ and $\omega$ are the connections on the gauge bundle $V$
  and the tangent bundle $TM$, respectively. 
Note that we can use $p=3$ in (\ref{stp}) to get a degree seven class upon integration
over the circle
$
\int_{S^1} ev^* ch_{4}(E)
$,
which gives 
 $
 \frac{1}{(2\pi)^3} \frac{1}{4!} C^3(LV) + {\rm decomposables}
 + {\rm non-gauge~factors}$.
 
 \vspace{3mm}
 In this case the fivebrane class can be described as follows. We have, for $n \geq 5$,
 the isomorphism 
 \(
 \int_{S^1} \circ~ ev^* : H^8\left( B{\rm String}(n); \Z\right) \to 
 H^7\left( LB{\rm String}(n); \Z  \right).
 \label{evh7}
 \)
%
%
%
  Since $H^7(B{\rm String}(n) ;\Z)= 0$ and $H^7({\rm String}(n); \Z)=\Z$,
  we have that the image in (\ref{evh7}) is $\Z$.   
    In terms of the space itself, the evaluation map and integration over the circle
    give 
    \(
    \int_{S^1} ev^* : H^8(M; \Z) \to H^7(LM; \Z).
    \)

%

 \vspace{3mm} 
 Next we consider the second case.
  In addition, replace the circle $S^1$ by a higher-dimensional space $Y$ so 
  that the loop space $LM=[S^1, M]$ is replaced by a higher degree generalization
 ${\rm Map}[Y,M]$. Then what replaces the evaluation map (\ref{ev1}) is
 $
 ev : Y \times {\rm Map}[Y,M] \to M
 $
 and (\ref{stp}) would then become  
 $
 \int_{Y} ev^* ch_{p+1}(E)
 $.
 The result will be a class on ${\rm Map}[Y,M]$ of degree $2p -{\rm dim}Y$. Obviously,
 when $p=1$ and ${\rm dim}Y=1$, we get back the String case. We will discuss further
 aspects of the general case in section \ref{6g}.
 
 \vspace{3mm}
 There are two special cases of interest, the first when $Y$ is a torus and the second 
 when $Y$ is a sphere. Let the dimension of $Y$ be $m$ and that of $X$ be $n$.
 Then the two cases are
 
  \begin{itemize}
 
\item  $Y=T^m$: This gives ${\rm Map}[T^m, M]=L^m M$, the higher 
iterated loop 
 spaces of $M$, i.e. $L^m M = L L \cdots L X^n$ ($m$ times). This is the 
 iterated loop space of $M$ which is obtained by looping on $X$ $m$
 times. The bundle replacing the loop bundle of the String case will be a bundle
 with structure group the {\it toroidal group} ${\rm Map}[T^m, G]=L^m G$.
  
 \item  $Y^m =S^m$: In this case the space to consider is ${\rm Map}[S^m; M]$,
 corresponding to the homotopy groups $\pi_m(M)$ of $M$. 
 Such spaces, at least for low $m$,  
 have been studied in connection to gauge theory in physics in \cite{Mic}.
\end{itemize}
 
 
 \vspace{3mm} 
In the physical situation under consideration, $Y$ can be taken to be the 
spatial part $W^5$ of the fivebrane worldvolume in spacetime of dimension
ten for the heterotic fivebrane and dimension eleven for the M5-brane.
Then, for $p=3$, the 
integration over the worldvolume yields a degree three class on the fifth loop
space $L^5 X^{10}$ (likewise for eleven dimensions). The 
 homotopy groups of the two groups, $G$ and $L^m G$ are related by
 $
 \pi_n(L^m G)=\pi_{n+m}(G) \oplus \pi_{n}(G)
 $,
so that, in particular, $\pi_3(L^5G)=\pi_8(G) \oplus \pi_3(G)$. For example, 7-connectedness
of $G$ would be the same as 2-connectedness of $L^5G$.

\subsection{Determinantal 6-gerbes}
  \label{6g}
  We can look at the relation to higher gerbes by looking at the structure {\it on} the
  worldvolume itself. Consider the case when the worldvolume has five compact 
  dimensions, i.e. time is considered to be $\mathbb{R}$.
 The Fivebrane structure can be
  described in terms of a 5-gerbe in general. In the special case of  
  an {\it index 5-gerbe}, this can be interpreted as the 
  degree five part of the Family's index theorem for Dirac operators
  \cite{Lo} \cite{BKTV}.
 
  \vspace{3mm}
  The construction of the index 1-gerbe is given in \cite{Lo} where also general 
  features of the higher gerbes are given. The construction for those is given
  in \cite{BKTV} which we follow below.  The degree two component of the 
  families index theorem is given by the first Chern class, or the curvature, 
  of a determinant line bundle.  There is an obstruction to realizing the component 
  of the families index theorem of degree higher than two as a curvature of some
  geometric object. 
  This, however, is automatic if we choose our spacetime $X$ to be 
  5-connected. so $X$ is homeomorphic to $S^{10}$, 
  and the explicit construction is given in \cite{Lo} 
  in this case.  The corresponding Deligne classes form a countable
  sets corresponding to different trivializations of the index bundle
  on the five-skeleton of the triangulation of $X$ and are labeled
  by $\bigoplus_{j=1}^2 H^{5-2j}(X;\Z)=H^1(X; \Z) \oplus H^3(X; \Z)$.

 \vspace{3mm}
   Assume $Y^m$ to be a compact oriented
  $C^{\infty}$-manifold of dimension $m$ over which we have a smooth  
  complex vector bundle $E$. Then
  the space of sections $\Gamma{E}$ of the bundle is expected to give rise to 
  a {\it determinental $5$-gerbe} 
  ${\rm Det} \Gamma(E)$, generalizing the determinant line bundle \cite{BKTV}.
  Thus, in the case of the fivebrane, if we take the spatial part 
  then the currently perceived wisdom leads us to a determinantal $5$-gerbe 
  and if we take the even-dimensional 
  spacetime then we are led to a determinantal $6$-gerbe. 
    
  \vspace{3mm}
  We now consider a family of fivebranes by considering the map to spacetime
  $q : W^5 \to X^{10}$ 
    of relative 
  dimension $5$ and a $C^{\infty}$ bundle $E$ on $W^5$. Then the characteristic
  class of the $5$-gerbe would be a class in 
  $H^{6}(X^{10}, C_{X^{10}}^{\infty *})=H^7(X^{10}, \Z)$. The class of the 
  determinantal $6$-gerbe in complex cohomology should be thought of
  as a $6$-fold delooping of the usual first Chern 
  (determinantal) class. 
  Still following \cite{BKTV}, the Real Riemann-Roch formula gives 
  \(
  C_1(q_* E)= \int_{W^5} \left[ ch(E) \wedge {\rm Td}(TW^5)\right]_{12} \in H^7(X^{10}, \C),
 \label{ind}
  \)
  where we are integrating the degree twelve part of the index formula over the 
  five-dimensional spatial part of the fivebrane worldvolume to get 
  a degree seven class.
  This involves the Dolbeault operator $\overline{\partial}$ over a complex
  envelope of $W^5$. In the smooth category, we just replace the Todd class
  Td in (\ref{ind}) with the $\widehat{A}$, the roof-genus of $TW^5$.

\section{The Gauge Algebra of Supergravity in $6k-1$ Dimensions}
\label{sec gauge alg}

Five-dimensional $SO(2)$ supergravity on a five-dimensional Spin manifold
$X$ is the theory obtained by coupling 
pure supergravity in five dimensions to an $SO(2)$ vector multiplet \cite{CN}\cite{C}. The
former is made of a metric $g$ on $X$ and a Rarita-Schwinger field $\psi$, which is a section 
of the spin bundle $SX$ coupled to the tangent bundle $TX$, 
$\psi \in \Gamma(SX\otimes TX)$. The latter contains an $SO(2)$-valued, hence abelian,
one-form $C_1$ with curvature two-form $G_2=dC_1$. The Lagrangian $\L_{(5)}$
will have

\begin{itemize}
\item a bosonic part $\L_{(5), {\rm bos}}$ for the bosonic fields $(g, G_2)$,
\item a fermionic part $\L_{(5), {\rm ferm}}$ for the fermionic field $\Psi$,
\item and an interaction part $\L_{(5), {\rm int}}$ for the terms that are mixed in the 
bosonic and fermionic fields.
\end{itemize}

In $\L_{(5)} =\L_{(5),{\rm bos}} + \L_{(5),{\rm ferm}} + 
\L_{(5),{\rm int}}$ we will consider only the 
bosonic part, given by the five-form
\(
\L_{(5),{\rm bos}}= R \, {*\oneone} -\frac{1}{2} G_2 \wedge *G_2 
-\frac{1}{6}G_2 \wedge G_2 \wedge C_1,
\label{5Lag}
\)
where $*$ is the Hodge duality operator on differential forms in five dimensions,
and $R$ is the scalar curvature of the metric $g$ of $X$.

\vspace{3mm}
Eleven-dimensional supergravity \cite{CJS} has some common features with five-dimensional
supergravity \cite{C} \cite{CN}, described above. The bosonic field content is the same, 
except that the potential $C_3$, replacing $C_1$, is now of degree three
so that the corresponding field strength $G_4$ is of degree four. The 
Hodge dual in eleven dimensions to $G_4$ is $G_7$. The bosonic part of the
Lagrangian is given by the eleven-form 
\(
\L_{(11),{\rm bos}}= R \, {*\oneone} -\frac{1}{2}  G_4 \wedge *G_4 
-\frac{1}{6}G_4 \wedge G_4 \wedge C_3\; .
\label{11Lag}
\)
From here on we treat both theories at the same time. We thus take $X$ to be a
$(6k-1)$-dimensional Spin manifold on which we define a supergravity with
Chern-Simons term built out of the potential $C_{2k-1}$, with corresponding
field strength $G_{2k}$. The value $k=1$ corresponds to the five-dimensional
case and the value $k=2$ to the eleven-dimensional case.

\vspace{3mm}
The equations of motion are obtained from the Lagrangian via the variational 
principle. The variation $\frac{\delta \L_{(6k-1),{\rm bos}}}{\delta C_{2k-1}}=0$ for 
$C_{2k-1}$ gives the corresponding equation of motion 
\(
d*G_{2k} + \frac{1}{2} G_{2k} \wedge G_{2k}=0 \; .
\label{d5 eom}
\)
We also have the Bianchi identity 
\(
dG_{2k}=0\;.
\label{5 Bian}
\)

\vspace{3mm}
The second order equation (\ref{d5 eom}) can be written in a first order form,
by first writing 
$
d\left( *G_{2k} + \frac{1}{2} C_{2k-1} \wedge G_{2k}\right) =0
$
so that 
\(
*G_{2k}= G_{4k-1} := dC_{4k-2}  -\frac{1}{2} C_{2k-1} \wedge G_{2k}\;,
\label{d5 1st}
\)
where $C_{4k-2}$ is the potential of $G_{4k-1}$,
the Hodge dual field strength to $G_{2k}$ in $6k-1$ dimensions.

\vspace{3mm}
The action $S_{\rm bos}=\int_{X} d{\rm vol}(X) \L_{\rm bos}$, and hence the
equations of motion, are invariant under the abelian gauge transformation
$\delta C_{2k-1}= d \lambda_{2k-2}$, where $\lambda_{2k-2}$ is a 
$(2k-2)$-form.  We can alternatively 
write the gauge parameter as $\Lambda_{2k-1}=d\lambda_{2k-2}$. 
In fact, the first order equation (\ref{d5 1st}) is invariant under the infinitesimal 
gauge transformations
\bea
\delta C_{2k-1}= \Lambda_{2k-1}, ~~~~~~~~~~~~~~~~~
\delta C_{4k-2}=\Lambda_{4k-2} -\frac{1}{2}\Lambda_{2k-1}
\wedge C_{2k-1},
\eea
where $\Lambda_{4k-2}$ is the $(4k-2)$-form gauge 
parameter satisfying $d\Lambda_{4k-2}=0$.
Applying two successive gauge transformations with different parameters $\Lambda_i$
and $\Lambda'_i$, $i=2k-1, 4k-2$, for the first and 
second one, respectively, and forming the the commutators gives
 \bea
 \left[ \delta_{\Lambda_{2k-1}}, \delta_{\Lambda'_{2k-1}}\right] 
 &=& \delta_{\Lambda''_{4k-2}}\; ,
 \nonumber\\
 \left[ \delta_{\Lambda_{2k-1}}, \delta_{\Lambda_{4k-2}}\right] &=& 0 \; ,
 \nonumber\\
 \left[ \delta_{\Lambda_{4k-2}}, \delta_{\Lambda'_{4k-2}}\right] &=& 0 \; ,
\label{5 comm lambda}
 \eea
with the new parameter $\Lambda''_{4k-2}=\Lambda_{2k-1} \wedge \Lambda'_{2k-1}$.
Note that the transformations are nonlinear, and this can be tracked back to the
presence of the Chern-Simons form in the Lagrangian (\ref{5Lag}).

\vspace{3mm}
We now introduce generators $v_{2k-1}$ and $v_{4k-2}$ for the $\Lambda_{2k-1}$ and 
$\Lambda_{4k-2}$ gauge transformations, respectively. On the generators, 
from the commutation relations (\ref{5 comm lambda}), we get the 
graded Lie algebra 
\bea
\left\{ v_{2k-1}, v_{2k-1} \right\} &=& - v_{4k-2} \;,
\nonumber\\
\left[ v_{2k-1}, v_{4k-2} \right] &=& 0 \;,
\nonumber\\
\left[ v_{4k-2}, v_{4k-2} \right] &=& 0 \;,
\label{5d v alg}
\eea
Note that we can use a graded commutator, which unifies a commutator and
an anticommutator, so that the above algebra (\ref{5d v alg}) becomes
\bea
\left[ v_{2k-1}, v_{2k-1} \right] &=& - v_{4k-2} \;,
\nonumber\\
\left[ v_{2k-1}, v_{4k-2} \right] &=& 0 \;,
\nonumber\\
\left[ v_{4k-2}, v_{4k-2} \right] &=& 0 \;,
\label{5d v grad alg}
\eea
where it is now understood that we are using graded commutators. The generators
satisfy the following properties
\begin{enumerate}
\item The generators $v_{2k-1}$ and $v_{4k-2}$ are constant: 
$dv_{2k-1}=0=dv_{4k-2}$.
\item The grading on the generators $v_{2k-1}$ and $v_{4k-2}$ follow 
that of the potentials
$A_{2k-1}$ and $A_{4k-2}$, respectively. Hence, $v_{2k-1}$ is 
odd and $v_{4k-2}$ is even. Thus,
$d(v_{2k-1} \alpha)= -v_{2k-1} d\alpha$ 
and $d(v_{4k-2} \alpha)= v_{4k-2} d\alpha$, for any 
$\alpha$.
\end{enumerate}
We will think of these ``generators" $v_i$ as elements of a graded Lie algebra,
where we will write $C_{2k-1}\otimes v_{2k-1}$, etc. instead of just 
$C_{2k-1}$ for the fields (see the discussion around equation (\ref{alg m})).

\vspace{3mm}
The field strengths can be combined into a total uniform degree field strength
$\mathcal{G}$ by writing
\(
\V = e^{C_{2k-1} \otimes \hspace{0.5mm}  v_{2k-1}} 
e^{C_{4k-2} \otimes \hspace{0.5mm}  v_{4k-2}},
\)
so that 
\bea
\G &=& dC_{2k-1} \otimes v_{2k-1} + 
(dC_{4k-2} - \frac{1}{2} C_{2k-1} \wedge dC_{2k-1}) \otimes v_{4k-2}
\nonumber\\
&=& G_{2k} \otimes v_{2k-1} + G_{4k-1} \otimes v_{4k-2}
\nonumber\\
&=& G_{2k}\otimes v_{2k-1} + *G_{2k} \otimes v_{4k-2}\; .
\eea
Note the analogy with usual (i.e. not higher-graded) nonabelian 
gauge theory. $\V$ is the analog 
of $g$ and $\G=d\V \V^{-1}$ is the analog of $dg g^{-1}$. 

\vspace{3mm}
The equation of motion for $C_{2k-1}$ ($=$ Bianchi identity for $C_{4k-2}$) 
and the 
Bianchi identity for $C_{2k-1}$ are obtained together from 
\(
d\G - \G \wedge \G = -d\V \wedge d\V^{-1} - d\V \V^{-1} \wedge d \V \V^{-1} =0\; .
\label{com 5}
\)
Indeed, using the commutators (\ref{5d v grad alg}), we have
\bea
\G \wedge \G &=& (G_{2k}\otimes v_{2k-1})\wedge( G_{2k}\otimes v_{2k-1}) + 
(G_{2k} \otimes v_{2k-1}) \wedge ( *G_{2k}\otimes v_{4k-2}) 
\nonumber\\
&&+(*G_{2k} \otimes v_{4k-2})\wedge (G_{2k} \otimes v_{4k-2})
\nonumber\\
&=& \frac{1}{2} [ G_{2k} \otimes v_{2k-1}, G_{2k} \otimes v_{2k-1}] 
+ [G_{2k} \otimes v_{2k-1}, *G_{2k} \otimes v_{4k-2}]
\nonumber\\
&=& \frac{1}{2} G_{2k} \wedge G_{2k} \otimes [v_{2k-1}, v_{2k-1}] 
- G_{2k} \wedge *G_{2k} \otimes [v_{2k-1}, v_{4k-2}]
\nonumber\\
&=& -\frac{1}{2} G_{2k} \wedge G_{2k}\otimes v_{4k-2} \;.
\label{gg 5}
\eea
Hence, (\ref{com 5}) follows from 
the equation of motion (\ref{d5 eom}) and the Bianchi 
identity (\ref{5 Bian}), which indeed correspond, respectively, 
 to the coefficient of $v_{2k-1}$ and $v_{4k-2}$ in the expression
 for $d\G$.
The case $k=2$, corresponding to eleven-dimensional supergravity, 
was derived in \cite{CJLP}.

\subsection{Models for the M-Theory Gauge Algebra}
In the previous section we have seen that $G_4$ and its dual 
can be written in terms of the total uniform degree field strength 
$\G$, the generators in which satisfy an algebra. It is natural to ask
about the nature of the generators and the graded structure in
which they result. In this section we provide a description in terms of 
homotopy (or higher-categorical) Lie  algebras:  $L_{\infty}$-algebras based on the 
constructions in \cite{SSS1}, and superalgebras 
corresponding to $(1|1)$ supertranslations.

\subsection{The gauge algebra as an $L_{\infty}$-algebra}

One connection to $L_{\infty}$-algebras is the appearance of
 higher form abelian Chern-Simons theory. 
 Recall that for $\gg$ any semisimple Lie algebra 
 and $\mu = \langle \cdot, [\cdot,\cdot]\rangle$
  the canonical 3-cocycle on it, we call
  $
    \gg_\mu
  $
  the corresponding (skeletal version of the)
 $\mathrm{string}_\mu(\gg)$  Lie 2-algebra. 
Similarly, for $\gg$ any semisimple Lie algebra and $\mu_7$
  the canonical 7-cocycle on it, we call
  $
    \gg_\mu
  $
  the corresponding (skeletal version of the) Fivebrane Lie 6-algebra
  \cite{SSS1} \cite{SSS3}.

\vspace{3mm}
\noindent{\bf Reminder on $L_\infty$-algebra valued diferential forms.}
Recall from \cite{SSS1} that for $\gg$ any $L_\infty$-algebra
 with $\mathrm{CE}(\gg)$ its Chevalley-Eilenberg differential graded 
 commutative algebra 
(DGCA), the space of GCA-morphisms $\mathrm{CE}(\gg) \to \Omega (X)$ is isomorphic
to the degree zero elements in the graded vector space
$
  \Omega^\bullet(X)\otimes \gg
 $,
where $\gg$ is in negative degree.
Flat $L_\infty$-algebra valued forms can be realized 
as graded tensor products 
$
  A \in \Omega^\bullet(Y)\otimes \gg
$
of forms with $L_\infty$-algebra elements with the special property that
\begin{itemize}
  \item
    $A$ is of total degree 0\,,
  \item
    $A$ satisfies a flatness constraint of the form
    \(
      d A + [A \wedge A] + [A\wedge A\wedge A ]
      + \cdots = 0
      \,,
    \)
    where $d$ and $\wedge$ are the operations in the 
    deRham complex and where $[\cdot, \cdots, \cdot]$
    are the $n$-ary brackets in the $L_\infty$-algebra.
    \end{itemize}
It is usually more convenient to shift $\gg$ by one into non-positive degree (hence with the usual Lie 1-algebra part in degree 0) and accordingly take $A$ to be of total degree 1.

\vspace{3mm}
Recall one of the central constructions of \cite{SSS1} involving 
the Weil algebra $\mathrm{W}(\gg)$.
For $\gg$ any $L_\infty$-algebra with degree $(n+1)$-cocycle
$\mu$ that is in transgression with an invariant polynomial $P$
\(
  \xymatrix{
    0&&
    P
    &&
    P
    \ar@{|->}[ll]
    \\
    \mu
    \ar@{|->}[u]_{d_{\mathrm{CE}(\gg)}}
    &&
    \mathrm{cs}
    \ar@{|->}[u]_{d_{\mathrm{W}(\gg)}}
    \ar@{|->}[ll]
    \\
    \mathrm{CE}(\gg)
    &&
    \mathrm{W}(\gg)
    \ar@{->>}[ll]
    &&
    \mathrm{W}(\gg)_{\mathrm{basic}}
    \ar@{_{(}->}[ll]
  }
\)
we can form the String-like extension
Lie $n$-algebra $\gg_\mu$ and the corresponding
Chern-Simons Lie $(n+1)$-algebra
$\mathrm{cs}_P(\gg)$ with the property that
$
  \mathrm{W}(\gg_\mu)
  \simeq
  \mathrm{CE}(\mathrm{cs}_P(\gg))
 $.
In \cite{SSS1} this construction was of interest for the case
that $\mu$ was a nontrivial cocycle on a semisimple Lie 1-algebra.
Another interesting case in which the construction works is when 
an invariant polynomial $P$ suspends to $0$, i.e. if it is in transgression with the 0-cocycle 
$\mu = 0$.

\vspace{3mm}
Notice that in particular all decomposable invariant polynomials
$P = P_1 \wedge P_2$, for $P_1$ and $P_2$ nontrivial and with transgression elements $\mathrm{cs}_i$,  $d_{\mathrm{W}(\gg)} \mathrm{cs}_i = P_i$ , suspend to 0, since for them we can choose the Chern-Simons element $\mathrm{cs} = \mathrm{cs}_1 \wedge P_2$,
which vanishes in $\mathrm{CE}(\gg)$ because $P_2$ does, by definition:
\(
  \xymatrix{
    0&&
    P_1 \wedge P_2
    &&
    P_1 \wedge P_2
    \ar@{|->}[ll]
    \\
    0
    \ar@{|->}[u]_{d_{\mathrm{CE}(\gg)}}
    &&
    \mathrm{cs}_1 \wedge P_2
    \ar@{|->}[u]_{d_{\mathrm{W}(\gg)}}
    \ar@{|->}[ll]
    \\
    \mathrm{CE}(\gg)
    &&
    \mathrm{W}(\gg)
    \ar@{->>}[ll]
    &&
    W(\gg)_{\mathrm{basic}}
    \ar@{_{(}->}[ll]
  }
  \,.
\)

\vspace{3mm}
\noindent{\bf Higher abelian Chern-Simons forms.}
A very simple but useful example are the decomposable invariant polynomials 
on shifted $\uu(1)$ in an even number of shifts: $b^{2k-2}\uu(1)$, for $k$ any 
positive integer. In this case 
\(
\mathrm{CE}(b^{2k-2}\uu(1)) =
(\bigwedge^\bullet (\underbrace{\langle c\rangle}_{2k-1}),
d= 0)
\)
and
\(
\mathrm{W}(b^{2k-2}\uu(1)) =
(\bigwedge^\bullet (
  \underbrace{\langle c\rangle}_{2k-1}
  \oplus
  \underbrace{\langle g\rangle}_{2k}
 ),
d c = g, dg = 0)
\,.
\)
The invariant polynomials are all the wedge powers
$g,\,,g\wedge g,\, g\wedge g \wedge g$
of the single indecomposable one $P := g$
\(
\mathrm{inv}(b^{2k-2}\uu(1)) =
(\bigwedge^\bullet (
  \underbrace{\langle P\rangle}_{2k}
 ),
d P  = 0)
\,.
\)
Notice that in this case the CE-algebra of the ``String-like extension'' $b^{2k-2}\uu(1)_{\mu = 0}$ is that of
$b^{2k-2}\uu(1) \oplus b^{4k-1}\uu(1)$:
\(
  \mathrm{CE}(b^{2k-2}\uu(1)_{\mu = 0})
  =
  \mathrm{CE}(b^{2k-2}\uu(1) \oplus
   b^{4k-3}\uu(1))
  \,.
\)

\vspace{3mm}
\noindent{\bf Abelian Chern-Simons $L_\infty$-algebras.}
Let $k \in \mathbb{N}$ be a positive integer. 
Then the Lie $(2k-2)$-algebra $b^{2k-2} \uu(1)$
has a decomposable degree $4k$ invariant polynomial
$P_{4k}$ which is the product of two copies of the standard
degree $2k$-polynomial. The corresponding Chern-Simons
Lie $(4k-1)$-algebra $\mathrm{cs}_{P_{4k}}(b^{2(k-1)}\uu(1))$
is given by the Chevalley-Eilenberg algebra of the form
\(
  \mathrm{CE}(\mathrm{cs}_{P_{4k}}(b^{2k-2)}\uu(1)))
  =
  \left(
    \bigwedge^\bullet(
        \langle
          c_{2k-1},
          g_{2k},
          c_{4k-2},
          g_{4k-1}
        \rangle
    ),
    d
  \right)
\)
where
\bea
  d c_{2k-1} &=& g_{2k}\; ,
\nonumber\\
  d c_{4k-2} &=& c_{2k-1} \wedge g_{2k} + g_{4k-1}\; ,
\nonumber\\
  d g_{2k} &=& 0\; ,
\nonumber\\
  d g_{4k-1} &=& g_{2k} \wedge g_{2k}
  \,.
\eea
This has a canonical morphism onto
\(
  \mathrm{CE}(b^{2k-2}\uu(1)\oplus b^{4k-3}\uu(1))
  =
  \left(
    \bigwedge^\bullet(
        \langle
          c_{2k-1},
          c_{4k-2}
        \rangle
    ),
    d = 0
  \right)
\)
with respect to which we can form the invariant or basic polynomials
\(
  \xymatrix{
    \mathrm{CE}(b^{2k-2}\uu(1)\oplus b^{4k-3}\uu(1))
   &
    \mathrm{CE}(\mathrm{cs}_{P_{4k}}(b^{2k-2)}\uu(1)))  
    \ar[l]_{\quad i^*}
    &
    ~\mathrm{basic}(i^*)
  \ar@{_{(}->}[l]
  }
  \,.
\)
This is the DGCA
  \(
    \mathrm{basic}(i^*)
    =
    \left(
      \bigwedge^\bullet(
         \underbrace{\langle g_{2k} \rangle}_{2k}
      \oplus
         \underbrace{\langle g_{4k-1} \rangle}_{4k-1}
      )
      ,
      \left(
        d g_{2k} = 0\,,
        d g_{4k-1} = g_{2k} \wedge g_{2k}
      \right)
    \right)
    \,.
  \)

We consider the $L_\infty$-algebra $\mathfrak{sa}$ which admits the above as its
Chevalley-Eilenberg algebra.
This $\mathfrak{sa}$ is a graded Lie algebra with generators $v_3$ and $v_6$ in degree 3 and 6, respectively, and with the graded Lie brackets being
\bea
  \left[ v_3,v_3\right]  &= & v_6 
\nonumber\\
  \left[ v_3, v_6 \right] & = & 0
\nonumber\\
  \left[ v_6,v_6 \right]  & =&  0
  \,.
  \label{alg m}
\eea
Flat differential form data with values in this $L_\infty$-algebra
is given by a degree 1-element
$
  A = G_{4} \otimes v_3 + G_7 \otimes v_6
  \in \Omega^\bullet(X)\otimes \gg
  $,
where
\begin{itemize}
  \item
    $G_{2k}$ a closed $2k$-form;
  \item
    $G_{4k-1}$ is a $4k-1$-form satisfying
    $
      d G_{4k-1} = G_{2k} \wedge G_{2k}
      $.
\end{itemize}
In total we have that a Cartan-Ehresmann connection with respect to $i^*$ is given by differential form data as follows:
\(
  \raisebox{70pt}{
  \xymatrix@C=8pt{
    \Omega^\bullet_{\mathrm{vert}}(Y)
    &&
    \mathrm{CE}(b^{2k-2}\uu(1)\oplus b^{4k-3}\uu(1))
    \ar[ll]_<<<{A_{\mathrm{vert}}}
    \\
    \\
    \Omega^\bullet(Y)
    \ar@{->>}[uu]
    &&
    \mathrm{CE}(\mathrm{cs}_{P_{4k}}(b^{2k-2)}\uu(1)))  
    \ar@{->>}[uu]^{i^*}
    \ar[ll]_<<<<<<{(A, F_A)}
    &
    \mbox{
      \begin{tabular}{l}
       $G_{2k} = d C_{2k-1}$
       \\ 
       $G_{4k-1} = dC_{4k-2} + C_{2k-1}\wedge G_{2k}$
     \end{tabular}
    }
    \\
    \\
    \Omega^\bullet(X)
    \ar@{^{(}->}[uu]
    &&
    \mathrm{basic}(i^*)
    \ar@{^{(}->}[uu]
    \ar[ll]
    &
    \mbox{
      \begin{tabular}{l}
       $d G_{2k} = 0$
       \\ 
       $d G_{4k-1} = G_{2k}\wedge G_{2k}$
     \end{tabular}
    }
  }
  }
  \,.
\)
This can be regarded as a certain Cartan-Ehresmann connection for 
the product of a line $(2k-1)$-bundle and a line $4(k-2)$-bundle

\vspace{3mm}
\noindent{\bf The situation for 11-dimensional supergravity.}
The local gauge connection data of 11-dimensional supergravity is
given by a 3-form $C_3$ with curvature 4-form
$G_4 = d C_3$ which can be captured in a duality-symmetric manner by
regarding $C_3$ as the data giving a flat connection
with values in
the abelian Chern-Simons Lie 6-algebra obtained by
setting $k=2$, subject to a self-duality constraint:
$
  A = G_4 \otimes v_3 + (* G_4)\otimes v_6
 $.
The flatness condition satisfied by this is then
equivalent to the equations of motion for $G_4$
\(
  (d A + [A \wedge A] = 0)
  \Leftrightarrow
  \left\lbrace
    \begin{array}{l}
       d G_4  = 0 \\
       d * G_4 = -\frac{1}{2}G_4 \wedge G_4
    \end{array}
  \right\rbrace
\)
Therefore, 
\begin{theorem}
The $C$-field and its dual in M-theory define
an $L_{\infty}$-algebra as their gauge algebra. 
\end{theorem}

\subsection{The gauge algebra as a Superalgebra}

We next interpret the M-theory gauge 
algebra in another novel way. The commutator of 
$v_{3}$ with $v_{6}$ and that of $v_{3}$ and $v_{6}$
are zero. Furthermore, the commutator of two $v_{3}$'s 
gives $v_{6}$, it is natural to suspect that each one of the two 
generators belongs to a different subspace in some grading.  
Indeed, these are the even and odd gradings, and we have

\begin{proposition}
The generators $v_{3}$ and $v_{6}$ form a Lie superalgebra 
of translations in $(1|1)$ dimensions.
\end{proposition}

\begin{remark}
This is analogous to the generator $\frac{\partial}{\partial \theta} + \theta 
\frac{\partial}{\partial x}$, where $x$ is an even coordinate and
$\theta$ is an odd coordinate.   
Thus we see an analog of the supersymmetric quantum mechanical relation
$\left\{ Q, Q \right\}={\mathcal H}$, where $Q$ is the supercharge and ${\mathcal H}$ 
is the Hamiltonian of the system. 
\end{remark}

\section{Duality-Symmetric Twists}
\label{duality sym twist}
In the twisted cohomology setting one can form uniform degree 
expressions for both the fields, e.g. the cohomology classes, 
and the twisted differential. In this section we consider twists of 
degrees higher than the familiar three. Also, given the discussion
of uniform degree fields in the previous section, it is natural to 
use these for more exotic twists.

\subsection{Degree seven twists}
\label{sec7}
The NS $H$-field in type II string theory serves as the 
twist in the K-theoretic classification of the RR fields.
This involves unifying the fields of all degrees into one total RR field. 
We investigate, based on \cite{7twist}, the case of heterotic string theory where
there are $E_8 \times E_8$ or ${\rm Spin}(32)/\Z_2$ gauge fields
in addition to the $H$-field. 

\vspace{2mm}
\noindent {\bf a. Rationally:} Considering the gauge field and its dual as a unified 
field, the equations of motion at the rational level contain a twisted 
differential with a novel degree seven twist.
Consider the case where the Yang-Mills group $G$, 
is broken down to an abelian
subgroup, thus making the curvature $F_2$ be simply $dA$. 
The manifolds $M^{10}$ are 
chosen such that this breaking via Wilson
lines (= line holonomy) is possible. The result of the variation of the action 
 $S=\int H_3 \wedge * H_3 + \int F_2 \wedge * F_2$
 with respect to $A$ gives $(d - H_7\wedge){\mathcal F}=0$, 
 where ${\mathcal F}=F_2 + * F_2$ is defined as the combined curvature,
 and $H_7$ is equal to $*H_3$ at the rational level, in analogy with the RR fields.
In this analysis, following \cite{7twist},
we used the ``Chapline-Manton coupling" $H_3=CS_3(A)$,
where $CS_3(A)$ is the Chern-Simons three-form for the connection
$A$, whose curvature is $F_2$. This gives a twisted differential $d_{H_7}=d-H_7\wedge$
which in nilpotent, i.e. squares to zero, $d_{H_7}^2=0$, since $H_7$ is closed. 

\vspace{2mm}
\noindent {\bf (i)} Let 
$v_n$ denote the $n$th generator of the complex oriented cobordism ring.
Consider the case $n=2$ and let 
$R = \R[[v_2, v_2^{-1}]]$ be a graded ring. The generator $v_n$ has dimension 
$2p^n-2$, so that at the prime $p=2$, $v_2$ has dimension 6. Let $d_{H_7}=d-v_2^{-1}H_7$ 
be the twisted de Rham differential of uniform degree one. Denote by 
$\Omega_{d_{H_7}}^i(M^{10}; \R)$ 
the space of $d_{H_7}$-closed
differential forms of total degree $i$ on $M^{10}$. The total curvature ${\mathcal F}$
is an element of degree two, $i=2$, in the above space of forms.
The equation $d_{H_7}{\mathcal F}=0$ defining the complex is just the
Bianchi identity and the equation of motion of the separate fields. 
Another possibility is to use the combination 
$\mathcal{F}=u_1^{-1} F_2 + u_2^{-1}F_8$,
where $u_1$ has degree two and $u_2$ has degree eight. 

\vspace{2mm}
\noindent {\bf (ii)} The second step is to ask whether the argument at the level of rational cohomology 
generalizes to some rational generalized cohomology theory.
A twisting of complex K-theory over $M$ is a principal $BU_{\otimes}$-bundle 
over $M$. From $BU_{\otimes} \equiv K(\Z, 2) \times BSU_{\otimes}$, 
the twisting is a pair $\tau=(\delta, \chi)$ consisting of a
determinantal twisting $\delta$, which is a $K(\Z,2)$-bundle over
$M$ and a higher twisting $\chi$, which is a $BSU_{\otimes}$-torsor.
Twistings are classified, up to isomorphism, by a pair of classes
$[\delta]\in H^3(M,\Z)$ and $[\chi]$ in the generalized cohomology group
$H^1(X,BSU_{\otimes})$. The former is  twisted 
K-theory, where the 
twist is given by the Dixmier-Douady 
(DD) class. The 
twistings of the 
rational K-theory of $X$ are classified, up to isomorphism, by \cite{Te} 
the group  $\prod_{n>1} H^{2n+1} (M;\Q)$.
This shows that, in addition to the usual $H^3\Q$-twisting, one
can in principle have twistings from $H^5\Q$ and $H^7\Q$ etc.
It is thus possible that a degree seven twist comes from complex K-theory.
However, it is not obvious how to isolate just the degree seven part from
the tower of all $n>3$ odd-dimensional twists.


\vspace{2mm}
\noindent {\bf b. Integrally:} The above generalizes the usual degree three 
twist that lifts to twisted K-theory and raises the natural question of whether at the integral 
level the abelianized gauge fields belong to a generalized cohomology theory.

\vspace{2mm}
\noindent {\bf (i)} The appearance of the higher degree generator connects nicely with the discussion in 
\cite{KS1}\cite{ KS2} \cite{KS3} on generalized cohomology in type II (and to some extent type I) string 
theories. Further the appearance of the $w_4=0$ condition in \cite{KS1}, interpreted in 
\cite{S4} in the context of F-theory, which is a condition in heterotic string theory, is another 
hint for the relevance of generalized cohomology in the heterotic theory.
Given the appearance of elliptic cohomology through $W_7=0$ and the $H_7$-twist, 
then the condition $W_7 + [H_7]=0$ is expected to make an appearance, which would
give rise to some notion of twisted structure in a similar way that the analogous condition 
$W_3 + [H_3]=0$ of \cite{FW} amounts to a twisted ${\rm Spin}^c$ structure.
This structure that we seek would be related to a twisted String structure, but is 
not quite the same but is implied by it, since
the String orientation condition implies $W_7=0$ via the action of the 
operation $Sq^3$. 
Further, we expect the modified condition to correspond to a differential $d_7$ in 
the AHSS of 
twisted generalized theories, possibly Morava K-theory and elliptic cohomology, since 
the `untwisted' differential is the first nontrivial differential there, in analogy to the 
`twist' for $d_3$ generated by $[H_3]$ in the K-theory AHSS.
 \begin{conjecture}
The cohomology class $W_7 + [H_7]$ corresponds
to a differential in twisted Morava K-theory and twisted Morava $E$-theory,
where $[H_7]$ acts as the twist. 
 \end{conjecture}
 This is a generalization of the statement that $W_3 + [H_3]$ 
 corresponds to the differential $d=Sq^3 + [H_3]\cup$ in 
 twisted K-theory.  
 Of course the construction of such twisted generalized 
 theories is not yet established. Nevertheless we note
 the following.
 
 \vspace{2mm}
\noindent  {\bf (ii)} In general, twists of a cohomology theory $E$ are classified by $BGL_1(E)$,
i.e. the twisted forms of $E^*(X)$ correspond to homotopy classes of maps 
$[X, BGL_1(E)]$. An equivalent way of saying this, which the more familiar one 
in the context of K-theory, is that the twists are classified by $B{\rm Aut} (E)$, 
where ${\rm Aut}(E)$ is the automorphism group of $E$. 
The homotopy groups of $BGL_1(E)$ are given as units of the ring
$E^0({\rm pt})$ in degree 1, and as $E^{k-1}({\rm pt})$ in degree $k>1$. For K-theory,
$K^0({\rm pt})=\Z$ gives $\pi_0 BGL_1(K)=\Z/2$,  $K^2({\rm pt})=\Z$ gives 
$\pi_3 BGL_1(K)=\Z$, which is detected by a map $K(\Z,3) \to BGL_1(K)$ giving the 
standard degree twist. In addition, there is $\pi_7 BGL_1(K)= \Z$.

\vspace{2mm}
\noindent  {\bf (iii)} Another possibility is the following. Twists of $TMF$ are classified by 
$BGL_1(TMF)$. This may have nontrivial homotopy in degree 7 coming from the homotopy 
in degree 6 of the connective theory $tmf$. 
\footnote{The homotopy groups of the spectrum $tmf$ and those of its periodic 
version $TMF$ are related
as $\pi_* (TMF)=\pi_* (tmf) \left[ (\Delta^{24})^{-1} \right]$, where $\Delta$ is the 
discriminant of elliptic curves.}

\subsection{Duality-symmetric twists in ten-dimensional string theory}
In type II string theory one
encounters the Ramond-Ramond (RR) fields
$
F=\sum_i F_{i}u^{-i}$,
where $u$ is the Bott generator and $i$ is the degree of the RR field, 
which is even for type IIA and odd for type IIB. This satisfies the 
equation
$
d_{H_3} F=0$,
where $d_{H_3}$ is the twisted differential, whose uniform degree expression
is
$
d_{H_3}=d + u^{-1} H_3\wedge
$.
Here $H_3$ is the Neveu-Schwarz (NS) 3-form. This is explained very well in \cite{Fr}. 

\vspace{3mm}
As we saw in section \ref{sec7},
in \cite{7twist} a degree seven twist was uncovered in heterotic and type I 
string theory, where the twist is given by the dual $H_7$ of the usual  
NS $H$-field $H_3$ in ten dimensions. The differential is of the form 
$
d_{H_7}= d + v^{-1} H_7 \wedge
$.
In this theory, for instance for $v=u^3$, one can form the 
uniform total degree one field strength \cite{Fr}
\(
H= u^{-1}H_3 + u^{-3}H_7\;,
\label{37}
\)
with corresponding potentials, or $B$-fields, of total degree zero
$
B=u^{-1}B_2 + u^{-3}B_6
$.

\vspace{3mm}
Given that the total field strengths are built of more than one component, 
we can ask whether the corresponding differential of uniform degree
might be built out of a twist that has more than one component. Consider a
candidate twisted de Rham differential with an 
expression of the form 
\(
d_H=d + u^{-1}H_3\wedge + u^{-3}H_7\wedge\;.
\label{dH}
\)
The square is $d_H^2$ contains the terms that are zero because $d_{H_3}$ 
and $d_{H_7}$ are differentials. In addition, there is the cross-term
$u^{-4}(H_3 \wedge H_7 + H_7 \wedge H_3)$, which is zero by antisymmetry of the
wedge product. Of course, another way to immediately see this is to write
$d_H$ as $d_{H_3}+ u^{-3}H_7 \wedge$ or as $d_{H_7}+ u^{-1}H_3 \wedge$.
Thus one can build a twisted graded de Rham complex out of such a differential.

\vspace{3mm}
\noindent {\bf Remarks.}
{\bf 1.} In fact, one can build a differential by adding to $d_H$ all expressions of the
form $u^{-i}H_{2i+1}\wedge$, i.e. 
\(
 d_H'=d + \sum_{i=0}^{\infty} u^{-i}H_{2i+1}\wedge\;.
\label{dh'}
\)

\noindent {\bf 2.} As differential forms, the $u$ are constant, i.e. $du=0$.
We can conceive of two modifications of this: First consider the generators 
appearing in front of the $H_{2i+1}$ to be independent. For example, in
\cite{7twist}, instead of (\ref{37}), we used the expression 
\footnote{
In section 
\ref{sec gauge alg} we used $v_i$ to indicate a generator of 
degree $i$, so to make a distinction we are using the notation
$v_{(i)}$ to indicate a generator of level $i$ in complex-oriented
generalized cohomology. We understand that the first notation
is more standard for the second notion, but since this is the only
occurrence of the higher Bott 
generators, then we hope it will not cause a confusion. 
}
\(
H=v_{(1)}^{-1}H_3 + v_{(2)}^{-1}H_7\;,
\label{v1v2}
\) 
where $v_{(1)}$ is still the Bott generator 
and $v_{(2)}$ is the generalization of that generator, i.e. identified as coming
from a complex-oriented generalized cohomology theory at the prime $p=2$.  
In this case it still holds that $v_{(1)}$ and $v_{(2)}$ are constants as differential forms.

\vspace{3mm}
From the above discussion the following immediately follows. 
\begin{proposition}
There is a twisted graded de Rham complex with differential 
$d + \sum_{i=1}^{\infty}v_{(i)}^{-1}H_{2i+1}\wedge$
, provided the differential forms $H_{2i+1}$ are 
closed. The coefficients $v_{(i)}$ are constant as differential forms
and can be taken to be either
dependent or independent.
\end{proposition}
\noindent 
One can associate analytic torsion  \cite{MWu}
to this type of twisted de Rham complex  
and a spectral sequence for the corresponding twisted 
cohomology \cite{LLW}.

\subsection{Duality-symmetric twists in eleven-dimensional M-theory}
In M-theory the situation is much more interesting. In this case we will see
that we can have twists of even degree and the generators are not independent
in the sense that they satisfy relations. Some aspects of this discussion have
been observed in \cite{S1} \cite{S2} \cite{S3}. 

\vspace{3mm}
In the low energy limit of M-theory, in addition to the metric and the 
gravitino,
 there is the $C$-field $C_3$ with field 
strength $G_4$.  We can build a differential with $G_4$ a twist as follows. 
The square of the expression
$
d_{G_4}=d + v_{3} G_4 \wedge
$
is
\footnote{In this section we have suppressed the tensor product between generators
and fields for ease of notation.}
\(
d_{G_4}^2=d^2 + d(v_{3} G_4 \wedge) +
 v_{3} G_4 \wedge d +  v_{3} G_4 \wedge v_{3} G_4 \wedge\;.
\label{dg4}
\)
On the right hand-side of (\ref{dg4}), the first term is always zero since
the bare $d$ is the de Rham differential. For the second term we need to 
decide whether $v_{3}$ is even or odd as a differential form. Since 
$G_4$ is even we see that we have to choose $v_{3}$ to be odd in 
order to cancel the third term. In addition, for the left-over from the second
term to be zero, $G_4$ has to be closed. The last term has no other term 
against which to cancel, so it has to be zero by itself. We need $v_{3}$ to be 
idempotent. This can be achieved either by the fact that the form degree
is odd or by the stronger condition that it squares to zero, i.e. that it is 
a Grassmann variable. The above discussion generalizes in an obvious 
way to the case when the coefficient has degree $2i-1$ and the field has degree
$2i$. Therefore we have
\begin{proposition}
The de Rham complex can be twisted by a differential of the form
$d + v_{2i-1}G_{2i}\wedge$ provided that $G_{2i}$ is closed and 
$v_{2i-1}$ is Grassmann algebra-valued. 
\end{proposition}

\vspace{3mm}
In M-theory one can consider the field dual to the $C$-field. This is a 
field strength $G_7$, which at the rational level is Hodge 
dual to $G_4$. We can use $G_7$ to twist the de Rham differential in 
the same way that $H_7$ did. Furthermore, in the same way as in 
(\ref{v1v2}) one could
form a duality-symmetric uniform degree field strength 
$
G=v_{3}^{-1}G_4 + v_{6}^{-1}G_7
$.
This expression can now be used to twist the de Rham differential, leading to 
\(
d_G=d+G\wedge=d+v_{3}^{-1}G_4\wedge + v_{6}^{-1}G_7\wedge\;.
\label{dgg}
\)
The conditions for (\ref{dgg}) to be a differential are given in the following.
\begin{proposition}
The de Rham complex can be twisted by the differential $d_G$ provided 
that either
\begin{enumerate}
\item $dG_7=0$ and $v_{3}^2=0$, or

\item $\{v_{3}, v_{3}\}=v_{6}$ and $dG_7=-\frac{1}{2}G_4\wedge G_4$.  
\end{enumerate}
Furthermore, this differential on $\G$ is equivalent to the equations of motion
and the Bianchi identity of the $C$-field. 
\label{propg}
\end{proposition}

\begin {remarks}
{\bf 1.} The first case would hold when there is no Chern-Simons term in the M-theory action.

\noindent {\bf 2.} The second case arises in M-theory and realizes the equation of motion for the $C$-field. 
While this is what appears in M-theory, mathematically we can have combination 
of even and odd fields of any degrees
\(
d+ v_{2m-1}^{-1}G_{2n} + v_{2m}^{-1}G_{2m+1}\;.
\) 
\noindent {\bf 3.}
The generators in proposition \ref{propg} have appeared in \cite{CJLP}
in the context of the M-theory gauge algebra (which is generalized in section
\ref{sec gauge alg}). What we have done above is relate
them to twisted cohomology. 
\end{remarks}

\section{M-brane Charges and Twisted Topological Modular Forms}

\subsection{Evidence for TMF}

\subsubsection{Construction of anomaly-free partition functions}

\vspace{3mm}
\noindent {\bf The $E$-theoretic partition function in type IIA.}

 The K-theoretic partition function encounters an anomaly \cite{DMW},
given by the seventh integral Stiefel-Whitney class $W_7$,
whose cancellation \cite{KS1}
is the orientation condition in elliptic cohomology 
for Spin manifolds, and in second integral Morava K-theory at $p=2$ 
$\widetilde{K}(2)$ 
(cf. \cite{Mo}) for oriented manifolds. 
The above class $W_7$ is the result of applying the 
Steenrod square operation $Sq^3$ on $w_4$, the 
fourth $\Z/2$ Stiefel-Whitney class or, equivalently, 
the 
result of applying the Bockstein operation $\beta=Sq^1$ 
on the degree six class $Sq^2 w_4$,
$
W_7=Sq^3 (w_4)= \beta  Sq^2 (w_4)
$,
by the Adem relation $Sq^3=Sq^1 Sq^2$.

\begin{theorem} [{\cite{KS1}}] 
\begin{enumerate} 
\item  A $10$-manifold $X$ is orientable with respect to $\tilde{K}(2)$ iff $W_7(X)=0$.
\item  The M-theory partition function is anomaly-free when constructed on $\tilde{K}(2)$-orientable
spaces.
\end{enumerate} 
\end{theorem}
Similar results hold also for Morava $E(2)$-theory.
In \cite{KS1} an elliptic refinement of the mod 2 index $j$ 
is obtained.  
Assuming 
that $X$ is orientable with respect to a real version 
$EO(2)$ of $E(2)$-theory, 
there is an $EO(2)$-orientation class
$[X]_{EO(2)}\in EO(2)_{10}(X)$. Now for $x\in E^0(X)$, the class
$x\overline{x}$ lifts canonically to $EO(2)^{0}(X)$, so 
$j(x)=\langle x\overline{x},[X]_{EO(2)}\rangle\in EO(2)_{10}$,
the right hand side being $EO(2)_{10}=\Z/2[v_{(1)}^{3}v_{(2)}^{-1}]$
by \cite{HK1}. The assumption on $EO(2)$-orientation is made precise:

\begin{theorem} [{\cite{KS1}}]
\begin{enumerate}
\item A spin manifold $X$ is orientable with respect to $EO(2)$ if 
and only if it satisfies $w_4(X)=0$, where $w_4$ is the fourth Stiefel-Whitney class.
\item When $w_4\neq
0$, this uncovers another anomaly to the existence of an elliptic cohomology
partition function.
\end{enumerate}
\end{theorem}

\begin{remark}
The class $w_4$ is the mod 2 reduction of the integral class
$\lambda=\frac{1}{2}p_1$, so that the vanishing of $\lambda$ implies
the vanishing of $w_4$, which in turn implies the vanishing of 
$W_7$. Therefore, the String orientation condition $\lambda=0$ 
is a necessary condition for the cancellation of the DMW anomaly.
\end{remark}

\noindent There is a ``character map'' $E\to K[[q]][q^{-1}] $ where 
$q$ is a parameter of dimension $0$, with $K[[q]]$ a
product of infinitely many copies of $K$, and 
the notation $[q^{-1}]$ signifies that $q$ is inverted. The map  
is determined by what happens on coefficients \cite{AHS}.  
\begin{theorem} [{\cite{KS1}}] {\it The refined partition function is a one-parameter family 
of theta functions.
At the prime $p=2$, $v_{(1)}^3v_{(2)}^{-1}$ is of dimension zero and serves as the expansion 
parameter $q$. }
\end{theorem}

\noindent{\bf The $E$-theoretic partition function in type IIB.}
The construction of the partition function in the IIB case is 
analogous to that of type IIA. Instead of $K^1(X)$, one starts with $E^1(X)$ where $E$
is a complex-oriented elliptic cohomology. The construction proceeds
precisely analogously as in the $K^1(X)$ case. However, the discussion of 
the phase is  delicate. First there is the pairing in $E^1(X)$: 
$E^1(X)\otimes E^1(X)\to E^2(X)\to E^{-8}=E^0$
where the second map is capping with the fundamental orientation class
in $E_{10}(X)$. To construct a $\theta$-function, a quadratic structure is needed,
which amounts to considering real elliptic cohomology: A product of an 
$x\in E^1(X)$ with itself can be given a real structure, which gives rise to an element of
$\omega(x)\in {\rm E}\R^{1+\alpha}X$, which when 
capped with the fundamental orientation class in ${\rm E}\R_{10}(X)$
gives an element in $E^{\alpha-9}$. This is a $\Z/2$-vector space generated by 
the classes $v_{(1)}^{3n-1}v_{(2)}^{2-n}\sigma^{-4}a^2,\; n\geq 1$\cite{HK1}.
Therefore, 
\begin{theorem} [\cite{KS2}] 
{\it There is a quadratic structure depending on one free parameter,
which leads to a precise IIB analog of the $\theta$-function constructed
for IIA using real elliptic cohomology in \cite{KS1}.}
\end{theorem}

\vspace{3mm}
\noindent {\bf TMF and the type IIB fields.}
A particularly convenient combination of the two degree-three fields in type IIB string theory
is ${G}_3=F_3-\tau H_3$, where $\tau$ is the parameter on the upper half plane.
This is a field with modular weight $-1$ since it transforms as
${G}_{3}^{\prime} ={G}_3\cdot(c\tau+d)^{-1}$ under ${\tau}^{\prime} =(a\tau+b)/(c\tau +d)$.
In $tmf$, a class of modular weight $k$ appears in $tmf^{2k}(X^{10})$. Therefore, 
\vspace{-1mm}
\begin{proposition} [{\cite{KS3}}]
{\it The fields of type IIB string theory as elements in 
tmf, satisfy ${\tilde{G}_3\in tmf^{-2}X^{10}}$.} 
\end{proposition}

This points to the $12$-dimensional picture: suppose, in the simplest 
possible physical scenario \cite{V} that $V^{12}=X^{10}\times E$  where $E$ is an 
elliptic curve, then $G_3 \times \mu \in tmf^{0}(V^{12})$ where $\mu\in tmf^2(E)$ 
is the generator given by orientation. It is consistent  that the class ends up in
dimension $0$ and no odd number shows up. Modular classes of weight $0$, however, must be 
in dimension $0$.
The mathematical interpretation of $\tau$ appears only when
we apply the forgetful map $E^k(X)\to K^k(X)[[q^{1/24}]][q^{-1/24}]$
with $q=\exp(2\pi i \tau)$. 
In fact, it is necessary to generalize to an
elliptic cohomology theory $E$ which is in general modular only
with respect to some subgroup $\Gamma\subset SL(2,\Z)$. 
For forms with such modularity, fractional powers of $q$ are needed: 
in the case of complex-oriented cohomology, one encounters $q^{1/24}$.
The map $E\to K[[q^{1/24}]][q^{-1/24}]$
whose induced map on coefficients (homotopy groups) makes the $k$-th 
homotopy group modular of weight $k/2$ is not the correct normalization 
to use because then $\tilde{E}^0(S^k)=E^{-k}({\rm pt})$ would
have modular weight $-k/2$ and not $0$. 
The correct normalization is given by composition with Adams 
operations (or alternatively with Ando operations \cite{And}) 
$\psi^{\eta}:K[[q^{1/24}]][q^{-1/24}]\to
K[[q^{1/24}]][q^{-1/24}]$, where $\eta$ is the Dedekind function
($\Delta^{1/24}$ where $\Delta$ is the discriminant form), which is a 
unit in $K[[q^{1/24}]][q^{-1/24}]$. For general $k$, multiplication
by $\eta^k$ is needed.

\vspace{3mm}
\noindent {\bf More on elliptic curves and F-theory.}
Roughly speaking, there is an elliptic cohomology theory for every elliptic 
curve. 
There is no universal elliptic curve over a commutative ring, so 
on what basis should one `favor' one elliptic curve over the
other, and hence one elliptic cohomology over the other? 
If one considers all elliptic curves at once then the corresponding
generalized cohomology theory is $tmf$, which is then, in a sense,
the universal `elliptic' cohomology theory. However, the price to pay
is that $tmf$ is {\it not} an elliptic cohomology theory and also 
not even-periodic. For more on this see e.g. \cite{Lu}.
In our identification of the F-theory elliptic curve with the elliptic 
curve in elliptic cohomology (after suitable reduction of coefficients),
we hence, along the lines of \cite{S4}, expect
 elliptic curves with a fixed modulus in F-theory to correspond 
to elliptic cohomology while ones with a modulus parameter varying 
in the base of the elliptic fibration, i.e. families,
 to correspond to $tmf$.   
It is then tempting to propose that $tmf$ sees all possible compactifications 
of F-theory on an elliptic curve, i.e. all admissible elliptic fibrations.

\vspace{3mm}
\noindent {\bf S-duality and twisted K-theory are not compatible.}
Type IIB string theory has a duality
symmetry, S-duality, which is analogous to 
electric/magnetic duality in
gauge theory. In the presence of $H_3$, 
the description of the RR fields of type IIB
using twisted K-theory is not immediately compatible
with S-duality. The origin of this is that the RR
fields are considered as elements of K-theory while
the S-dual field $H_3$ is taken to be a cohomology
class, leading to the breaking of the symmetry. 
 Furthermore, type IIB string theory has a five-form
in place of the four-form in type IIA and in M-theory. 
In \cite{DMW} the apparent puzzle about the incompatibility of 
twisted $K$-theory and S-duality
is raised. A definite statement is proved in \cite{KS2}.
The condition for anomaly cancellation
for $F_3$ is  $(Sq^3 + H_3)\cup F_3=0$, which is not
invariant under the full $SL(2,\Z)$ group. The direct $SL(2,\Z)$-invariant 
extension of the above equation is \cite{DMW}
$F_3 \cup H_3 + \beta Sq^2 (F_3 + H_3)=0$.
One immediate question is that of justification (and interpretation)
of the nonlinear term $\beta Sq^2H_3=H_3\cup H_3$.
The point is to exhibit this as a differential, or obstruction
for the cohomological pair $(H_3, F_3)$ to lift to the theory.

\vspace{3mm}
The usual requirement that twisted $K$-theory 
be a {\em module} over $K$-theory, which forces the `structure group' 
of the bundle of $K$-theories in question to be the multiplicative 
infinite loop space $GL_1(K)$ of $K$-theory, violates the condition.
One can then ask for some further generalized twisting, where, for 
a particular $H_3$, the choice of allowable $F_3$'s would not form a vector space,
i.e. whether one could consider a form of twisted $K$-theory
which is not a module cohomology theory over ordinary $K$-theory.
In \cite{KS2} this is shown not to exist if the twisting space
is $K(\Z,3)$.  
The {\em classifying space}, i.e. a topological space $\mathcal{B}$
such that the affine-twisted $K^1$-group would be classified by homotopy
classes of maps $X\to \mathcal{B}$, cannot exist because the `group cohomology'
$H^2(\C P^{\infty},BU)$ vanishes. The above equation 
cannot occur as first Postnikov 
invariant, so 
\begin{theorem} [{\cite{KS2}}] {\it  $K(\Z,3)$-twisted K-theory is not compatible with S-duality in
type IIB string theory.}
\end{theorem}

\noindent{\it 4. Re-interpreting the twist:}
Due to the above, one has to seek a solution in the realm of higher generalized 
cohomology theories. However, a question arises whether or not to leave the twist. 
The twist introduces an intrinsic non-commutativity 
which seems to prevent further delooping
of the theory into a modular second cohomology group, giving an indication for a  
 an {\em untwisted} generalized cohomology theory.
 
 \vspace{3mm}
 There is a map $K(\Z,3)\to TMF$ coming from the String orientation. 
 A twist of K-theory of the form $X \to K(\Z,3)$ then gives rise to an 
 element of $TMF(X)$ by composition $X \to K(\Z,3) \to TMF$, 
in fact defining ``elliptic line bundles''  (cf. \cite{D}). Let us explain this.  
The classifying space $B{\rm GL}_1 K$ for elementary twistings of complex $K$-theory
splits, as an infinite loop space, as a product of two factors $A \times B$. The first 
factor is a
$K(\Z,3)$ bundle over $K(\Z/2,1)$ which splits as a space but has nontrivial
infinite loop structure classified by $Sq^3 \in H^3(\H(\Z/2);\Z)$.
There is a natural infinite loop map $B \to TMF$ from
$B$ to the representing space for topological modular forms, and
so by projecting through $B$ a map $B{\rm GL}_1 K \to TMF$.  In particular
an elementary twisting of $K$-theory for $X$ determines a $TMF$-class on $X$
\cite{D}.
The geometric interpretation of these
$TMF$ classes is simplified if restricted to those
classes coming from twistings involving only the $K(\Z,3)$ factor of $B$.
Such a twisting is determined by a map $X \to K(\Z,3)$ or equivalently
by a $BS^1$ bundle on $X$.  This bundle can be thought of as a stack
\footnote{These are yet another incarnation of the line bundle gerbes
mentioned earlier, for instance in the paragraphs just before proposition
\ref{AJ prop}. This stack incarnation 
is precisely the ``gerbe" in the original sense of the word
(a locally non-empty and transitive stack).
}
 locally
isomorphic to the sheaf of line bundles on $X$ and hence as a
1-dimensional 2-vector bundle on $X$.  In this sense the
$TMF$ classes coming from $K$-theory twistings can be viewed 
as 1-dimensional elliptic
elements and twisted $K$-theory as $K$-theory with coefficients in
this ``elliptic line bundle'' \cite{D}.

\vspace{3mm}
Therefore, what looks like twisting
to the eyes of $K$-theory, untwists and becomes merely a multiplication
by a suitable element in $TMF$ or any suitable form of elliptic cohomology.
Thus,
\begin{observation} [{\cite{KS2}}] 
{\it  If both $F_3$ and $H_3$ are viewed as elements of elliptic cohomology,
i.e. symmetrically, and the twisting is replaced by multiplication then the 
S-duality puzzle is solved.}  
\label{S cov} 
\end{observation}

\subsection{Review of D-brane charges and twisted K-theory}
\vspace{3mm}
The analogy with the more familiar case of the NS field 
$H_3$ is as follows. The cohomology class $[H_3]$ appears in the definition
of a twisted ${\rm Spin}^c$-structure
\(
W_3 + [H_3]=0,
\label{spinc}
\)
a condition for consistent wrapping of D-branes around cycles in ten-dimensional 
spacetime
 \cite{FW},
 where $W_3$ is the third integral Stiefel-Whitney class of the normal 
 bundle, the vanishing of which
 allows a ${\rm Spin}^c$-structure.
In the presence of the NS B-field, 
or its
field strength $H_3$, the relevant K-theory is 
twisted K-theory, as was shown in \cite{Wi1} \cite{FW} \cite{Kap} by analysis of
worldsheet anomalies for the case the NS field  
$[H_3] \in H^3(X, \mathbb Z)$ is a torsion class, and in \cite{BM} for the
nontorsion case. 
Twisted K-theory has been studied for some time \cite{dk}  \cite{Ros}.  
More geometric flavors were given in \cite{BCMMS}. Recently, the theory
was fully developed by Atiyah and Segal \cite{AS1} \cite{AS2}. 
It is a further result that $[H_3]$ acts as a determinantal
  $K(\Z, 2)$-twist for complex K-theory. 
The left hand side of the expression (\ref{spinc}) is in fact the first differential $d_3$ 
in the Atiyah-Hirzebruch spectral sequence for twisted K-theory -- see  
\cite{BCMMS} \cite{AS1} \cite{AS2}. 
Then,
a {\it twisted D-brane} in a
$B$-field $(X,H_3)$ is a triple $(W,E,\iota)$, where
$\iota:W\hookrightarrow X$ is a closed, embedded oriented submanifold
with $\iota^*H=W_3(\W)$, and $E\in K^0(\W)$ (see \cite{BMRS}).

\vspace{3mm}
A D-brane wrapping a homology cycle is inconsistent if it suffers from an anomaly, and is sometimes inconsistent if the homology 
cycle cannot be represented by any nonsingular submanifold. This is detected by the first Milnor primitive cohomology operation
$Q_1=\beta P_3^1$ where $\beta$ is the Bockstein and $P_3^1$ is the Steenrod power operation,
both at the prime $p=3$. In contrast, the twisted Spin${}^c$ condition is at $p=2$.
In fact, we have
\begin{observation} [{\cite{ES2}}]
{\it The twisted Spin${}^c$ condition is not sufficient.}
\end{observation}

\vspace{3mm}
 $D$-brane charges are
classified by the twisted K-group $K^*_{H}(X)$.
A rigorous formulation of such $D$-brane charges requires 
a Thom isomorphism and a push-forward map. Indeed
the Thom isomorphism and push-forward in twisted K-theory are 
established in \cite{CW}:
Corresponding to the map $\iota: \W \hookrightarrow X$ there is 
\begin{enumerate}
\item Push-forward map: $\iota_{!} : K^{\bullet}_{\iota^* \sigma+ W_3(\iota)}(\W)
\to K^{\bullet}_{\sigma}(X)$.

\item Thom isomorphism: $K^{\bullet}(\W) \cong K^{\bullet}_{\iota^* \sigma + W_3(\iota)}
(\W)$.
\end{enumerate}
With the use of the Riemann-Roch formula and index theorem in twisted
K-theory it is now established that the RR charges in the presence of an $H$-field
are indeed classified by twisted K-theory \cite{CW2}.

\vspace{3mm}
Now applying the push-forward map for $\iota$ in 
twisted K-theory,
one can associate a canonical element in $K_H(X)$, 
 the desired $D$-brane charge of the underlying $D$-brane 
 \cite{CW}
\(
\iota_!: K({\mathcal{W}}) \cong K_{\iota^*H + W_3(\W)} (\W) \longrightarrow
K^*_{H} (X),
\)
Hence, 
\begin{definition}
For any $D$-brane wrapping $W$ determined by an element
$E \in K(\W)$, the charge is
\(
\iota_!(E) \in K^*_{H}(X).
\)
\end{definition}

\subsection{The M-brane charges and twisted TMF}
\label{M ch tmf}
In the case of M-theory, the object
carrying charges with respect to $G_4$ is the M5-brane and we will study
conditions for consistent wrapping of such branes on cycles in 
eleven-dimensional spacetime. Hence, the candidate object to 
carry charges with respect to TMF is the M-theory fivebrane.
In this section we will provide a point of view on the 
interpretation of Witten's quantization condition (\ref{shift}) 
for $G_4$, which will 
give the context within which we describe M-brane charges
in the following section.

\vspace{3mm}
 In the case of the string, the 
target spacetime is assumed to be Spin, i.e. $w_2(X^{10})=0$. Then this
also implies that $X^{10}$ is certainly Spin$^c$. Then the requirement 
that the brane's worldvolume be Spin$^c$ is
equivalent to requiring the normal bundle to the D-brane to be Spin$^c$.
On the other hand, for the M-theory case, Witten's flux quantization is 
obtained from the embedding of the membrane in spacetime $Y^{11}$.
Taking $Y^{11}$ be Spin$^c$ will not be enough this time. We will see
below how the twisted String condition \cite{Wa} \cite{SSS3} 
in the three bundles: tangent bundle to the worldvolume, 
the normal bundle, and the tangent bundle of the target will be
related.  
\vspace{3mm}
\begin{remarks}
\label{rem 4}
{\bf The Interpretation of Witten's quantization.}

\noindent {\bf 1.}
 The equation (\ref{shift}) makes sense only if $\lambda$ is divisible by
two. It means that it is not $[G_4/2\pi]$ but $[G_4/\pi]$ that is well defined as 
an integral cohomology class and that this class is congruent to $\lambda$ modulo
two \cite{Flux}.  

\noindent {\bf 2.}  A weaker condition than (\ref{shift}) can be obtained by multiplying 
by two, provided there is no 2-torsion,
\(
2 [G_4] + \lambda = 2 a. 
\label{times2}
\)
Condition (\ref{times2}) gives condition (\ref{shift}) if there is no
2-torsion and once $\lambda$ is divisible by two.

 \noindent {\bf 3.} We rewrite condition (\ref{times2}) in turn in a suggestive way as
$
\lambda - 2\left( [G_4] - a \right)=0
$,
so that we identify $\alpha:=2([G_4] - a)$ as the twist of the String structure
\cite{Wa} \cite{SSS3}.
 
 \noindent {\bf 4.} Alternatively, we can work not with twisted String structure
 but rather with what was called twisted ${\mathcal F}^{\langle 4 \rangle}$-structure
 in \cite{SSS3} to account for the factor of 2 dividing $\lambda$. 
 There is no canonical description
 of ${\mathcal F}^{\langle 4 \rangle}$ yet except through $BO \langle 8 \rangle$.
\end{remarks}

\vspace{3mm}
We can consider 3 different cases: 
\begin{itemize}
\item the case when $a=0$ so that the 
$E_8$ bundle is trivial and the twist is provided by $[G_4]$,
\item the case $[G_4]=0$ so that the flux is $G_4=dC_3$ and the 
twist is provided by the $E_8$ class $a$,
\item the general case, where the twist is provided by $\alpha$. 
\end{itemize}

\vspace{3mm}
To make a comparison, 
let us briefly recall the model of \cite{DFM}.
The field strength in M-theory is geometrically described as a shifted
differential character \cite{DFM} in the sense of \cite{HS}. A shifted differential character is 
the equivalence class of a differential cocycle which trivializes a specific differential 
$5$-cocycle related to the integral Stiefel-Whitney class $W_5(Y)$. The 
Stiefel-Whitney class $w_4(Y) \in H^4(Y; \Z_2)$ defines a differential cohomology 
class $\check{w}_4$ via the inclusion 
$
H^4(Y;\Z_2) \hookrightarrow H^4(Y; \R/\Z ) \hookrightarrow H^4(Y;\Z )
$.
On a Spin manifold, $W_5(Y)=0$ is satisfied since $\lambda$ is an integral lift of 
$w_4(Y)$. In this case, the differential cohomology class $\check{w}_4$  can 
be lifted to a differential cocycle by defining 
$
\check{W}_5(Y)=\left(0, \frac{1}{2} \lambda, 0  \right)
\in \check{Z}^5(Y) \subset C^5(Y;\Z) \times C^4(Y;\Z) \times \Omega^5(Y)
$,
and the $C$-field can be defined as the differential cochain
$
\check{C} = (\overline{a}, h, \omega) \in C^4(Y;\Z) \times C^3(Y;\Z) \times \Omega^4(Y)
$
trivializing $\check{W}_5$, $\delta \check{C}=\check{W}_5$, i.e. in components 
\cite{DFM}
\(
\delta \overline{a} =0,
\qquad
\delta h = \omega - \overline{a}_{\R} + \frac{1}{2} \lambda,
\qquad
d \omega =0\;.
\)
It was proposed in in \cite{DFM} that $G_4$ lives in 
\(
\check{H}_{\frac{1}{2}\lambda}^4(Y^{11}),
\label{check H}
\)
 the space of 
shifted characters on $Y^{11}$ with 
shift $\frac{1}{2}\lambda$, and similarly on the fivebrane worldvolume. 

\begin{remark} In defining $G_4$ to live in (\ref{check H}), the authors
of \cite{DFM} are taking the point of view that the class $\frac{1}{2}\lambda$ 
acts as a twist for the differential character. The point of view we would like
to take here is that $\frac{1}{2}\lambda$ is what is being twisted, and hence
plays a more central role.  After all, natural structures, i.e. ones related to the
tangent bundle, should be in a sense more fundamental for describing 
structures on a manifold that are extra or auxiliary structures such as 
bundles not related to the tangent bundle. 
The $E_8$ gauge theory can then be seen as
responsible for the twist of $\frac{1}{2} \lambda$. This is the point of view also
adopted in \cite{SSS3}, building on the definition in \cite{Wa}. Therefore
\(
\left\{ \frac{1}{2}\lambda-{\rm shifted~} E_8 {~\rm structure}  \right\}
~~\Longrightarrow ~~
\left\{ E_8-{\rm twisted~String~ structure}  \right\}\; .
\)
\end{remark}

\subsection{The M5-brane charge}
There are two definitions for the the M5-brane partition function \cite{Effective} \cite{Among}
\cite{HS}, and hence for the M5-brane charge \cite{DFM}. One is 
intrinsic and uses the theory on 
the worldvolume. The other is extrinsic and uses anomaly inflow and hence
the embedding in eleven dimensions.
We will use the extrinsic definition for the M5-brane charge as this is the one that
uses the Thom isomorphism and the push-forward (in the appropriate theory). 
Some aspects of the intrinsic approach were used in \cite{SSS3} to relate
twisted Fivebrane structures to the worldvolume theory of the fivebrane.

\vspace{3mm}
Consider the embedding $\iota : {\mathcal W} \hookrightarrow Y$
 of the fivebrane with six-dimensional worldvolume ${\mathcal W}$ into eleven-dimensional
 spacetime $Y$. Consider the ten-dimensional 
 unit sphere bundle $\pi : V \to {\mathcal W}$ of ${\mathcal W}$ with fiber $S^4$ associated to the
 normal bundle $N \to {\mathcal W}$ of the embedding $\iota$. There is a corresponding 
 11-manifold $Y_r$ with boundary $V$ obtained by removing the disk bundle of
 radius $r$, $Y_r =Y- D_r(N)$ \cite{DFM}. 
 Corresponding to the sphere bundle $V$ there is the Gysin sequence. Using that
 the normal bundle has vanishing Euler class $e(N)=0$, 
 one can deduce that \cite{DFM}
 \bea
 H^3(V;\Z) &\cong& H^3({\mathcal W};\Z)
 \nonumber\\
 H^4(V;\Z) &\cong& H^4({\mathcal W};\Z) \oplus \Z~~~({\rm noncanonically}).
 \eea
 so that
 $
 H^3(V , U(1)) \cong H^3({\mathcal W}, U(1))
 $.
 This relates $G_4$ on the worldvolume to that on the normal bundle, i.e.
 the tangential components to the transverse components.

\vspace{3mm}
The M5-brane is magnetically charged under the $C$-field, i.e. the former 
acts as a source for the latter. The charge is then measured by the value of the
integral of $G_4$ over the linking sphere $S^4$ of ${\mathcal W}$ in $Y$
\(
Q_{M5}= \int_{S^4} G_4 =k \in \Z \;.
\)
Since the total Pontrjagin class of $S^4$ is 1, then 
$Q_{M_5}$ is equal to the instanton number of an $E_8$ instanton of the 
$E_8$ gauge theory with four-class $a$:
$
\int_{S^4} a=k \in \Z$.
The description in terms of instantons is further given in \cite{ES}. Note that
the charge defined this way highlights the role of $E_8$. If the $E_8$ bundle 
is trivial then the charge is zero. 

\vspace{3mm}
The anomaly inflow cancellation allows for a partition function of the M5-brane 
to be well-defined. In order for this to be made precise, the $C$-fields on ${\mathcal W}$ 
must be related to the $C$-field on $V$.  
The gauge equivalence classes of $C$-fields on $W$ is the shifted differential 
character $[\overline{C}] \in \check{H}_{\frac{1}{2}\lambda}^4(W)$, and similarly
for $W$. This then requires the existence of the map
\(
i : \check{H}_{\frac{1}{2}\lambda_W}^4(W) \to 
\check{H}_{\frac{1}{2}\lambda_X}^4(X).
\label{push check}
\)
In \cite{DFM}, a $C$-field $\check{C}_0$ on $X$ was {\it chosen} such that
via the map $i$ the $C$-fields are related as
\(
i [\overline{C}]=[\check{C}_0] + \pi^*[\overline{C}].
\label{C relations}
\)

\vspace{3mm}
We now give our main proposal.

\begin{proposal}
\begin{enumerate}
\item As we mentioned earlier, the shifted differential character 
gives the minor role to $\frac{1}{2} \lambda$. This obscures an interpretation 
in terms of generalized cohomology. We would like to replace 
(\ref{push check}) with
\(
\rho : \left(\frac{1}{2} \lambda + \alpha \right)({\mathcal W}) 
~\longrightarrow ~ \left(\frac{1}{2} \lambda + \alpha \right)(V) \; .
\)
\item
In order to properly define the M5-brane charge, what 
we need is then a Thom isomorphism and a push-forward
map in the appropriate theory. Our proposal is that the desired theory
is TMF and hence we apply the Thom isomorphism and the 
push-forward in TMF to obtain the M5-brane charge.  This way,
the analog of (\ref{C relations}) will be canonical. 
\end{enumerate}
\end{proposal}

\vspace{3mm}
If $X$ is a space, then the twisted forms of  $K^*(X)$ correspond to 
homotopy classes of maps $[X,BGL_1 K]$. The third homotopy group 
of the parametrizing space is $\pi_0 BGL_1 K = \Z$  since $K^2(pt) = \Z$.
This is the determinantal twist in K-theory an example of which being the 
NS $H$-field in string theory. 
Twists of TMF are classified by $BGL_1 TMF$ and there is a corresponding
map $K(\Z,4) \to BGL_1 TMF$. 

\vspace{3mm}
The proof of the following theorem is explained by Matthew Ando.
\begin{theorem} A class $\alpha \in H^4(X;\Z)$ 
gives rise to a twist $tmf_{\alpha}^* (X)$ of $tmf(X)$. 
Moreover, if $V$ is a (virtual)
spin vector bundle over $X$ with half-Pontrjagin class $\lambda$, 
then $tmf_{\lambda}^*(X) \cong tmf^*(X^V)$
as modules over $tmf^*(X)$.
 \end{theorem}

Then, armed with a Thom isomorphism and a pushforward map
(see Ando's contribution to these proceedings \cite{Atwist}),
the main application is

\begin{deftheorem} 
Given an embedding $\iota : {\mathcal W} \hookrightarrow Y$, 
the charge of the M5-brane is given by
\(
\iota_{!} (E) \in TMF_{\alpha}^*(Y) \; 
\)
\end{deftheorem}

\begin{remarks} We consider the Witten quantization condition 
for the worldvolume, normal bundle, and target for both the M2-brane
and the M5-brane.

\begin{enumerate}
\item {\it M2-brane}: The condition $\frac{1}{2}\lambda (M3) 
+ [G_4]|_{M^3} - a_{E_8}|_{M^3}=0$ is satisfied on the M2-brane worldvolume,
by dimension reasons. Therefore, the condition $\frac{1}{2} \lambda (N) 
+ [G_4]|_N - a_{E_8}|_N=0$ on the normal bundle $N$ is equivalent to the
same condition on the target. This is indeed what enters in Witten's
derivation of the quantization condition (\ref{shift}).

\item {\it M5-brane}: Assuming as above that the condition 
holds for the normal bundle to the fivebrane, then this implies that the
condition on the worldvolume is equivalent to that on the target. So we have
$$
\frac{1}{2} \lambda ({\mathcal W}) 
+ [G_4]|_{\mathcal W} - a_{E_8}|_{\mathcal W}=0
\hspace{1cm} \Longrightarrow
 \hspace{1cm}
\frac{1}{2} \lambda (Y) 
+ [G_4]|_Y - a_{E_8}|_Y=0\; .
$$
\end{enumerate}
\end{remarks}
 
\vspace{3mm}
\noindent {\bf Discussion and further evidence for the interpretation of 
M5-brane charge.}

\vspace{1mm}
{\bf 1.} {\bf Two-gerbes.}
The system of M2-branes ending on M5-branes is 
an M-theoretic realization of the system of strings ending 
on D-branes. In the latter, there is a gauge field, or connection 
one form $A$, on the string boundary $\partial \Sigma \subset Q$.
In addition, there is the $B$-field on $\Sigma$, which is 
acts as a twist for the Chan-Paton bundle, whose connection is
$A$ and curvature is $F$. Now in the case of M-theory, the
membrane boundary $\partial {\rm M}2 \subset {\rm M}5$ (written schematically)
has a degree two potential, essentially a $B$-field, 
which represents a gerbe made nonanbelian by the presence
of the pullback of the $C$-field. This is equivalent to a 2-gerbe 
system. 

{\bf 2.} {\bf Loop variables for the membrane.} 
The system of multiple M5-branes, generalizing the
system of $n$ D-branes leading to $U(n)$ nonabelian 
gauge symmetry, can be described by twisted $\widetilde{\Omega} G$-
gerbes, i.e. twisted gerbes for the universal central extension of the
based loop group $\Omega G$, 
where $G$ is any of the Lie groups ${\rm Spin}(n)$, $n \geq 7$, 
$E_6$, $E_7$, $E_8$, $F_4$, and $G_2$ \cite{AJ}. 
Indeed, it has been shown explicitly in \cite{G3} 
that classical membrane fields are loops. 

{\bf 3.} {\bf Fivebrane in loop space.}
The generalization of the abelian system of fields, called tensor multiplet, 
to the nonabelian case leading to the nonabelian tensor multiplet,
which appears in the worldvolume theory of the M5-brane,
requires loop space variables and a formulation in loop space 
\cite{G1} \cite{G2}. Given the general principle that degree $n$ phenomena 
in a space are captured by degree $n-1$ phenomena on its
loop space, the situation for $n=3$ suggests a relation between 
(twisted) K-theroy in loop space to be related to (twisted) 
$TMF$-cohomology of the space.

{\bf 4.} {\bf Relating TMF in M-theory to twisted K-theory in string theory.}
As we have reviewed it is known that twisted K-theory 
classifies D-branes and their charges in the presence 
of the NS $B$-field on a ten-dimensional space $X^{10}$. 
We have also seen how M-branes and their
charges should take values in TMF on an eleven-dimensional
space $Y^{11}$. Given the relation between M-theory and type IIA
string theory, the situation when $Y^{11}$ is a (possibly trivial) 
circle bundle over $X^{10}$, there should be a relation between 
the TMF description and the twisted K-theory description, in the
sense that (possibly $S^1$-equivariant twisted )TMF of $X^{10} \times S^1$ 
should give rise to twisted K-theory of $X^{10}$. 
Current discussions with Matthew Ando 
and with Christopher Douglas 
suggest schematically the following
\begin{conjecture}
There is a map $\left( \Omega tmf \times S^1\right) /S^1 \to k/*$.
\end{conjecture}

\noindent {\bf 5.} {\bf Capturing the fields of degree $4k$ in M-theory.}
Twisted K-theory, under the twisted Chern character map that lands in 
rational cohomology, leads to 
differential forms of all even degree up to the dimension of the manifold. 
These forms are the components of the RR field in the classical 
supergravity approximation. In eleven dimensions, then, one 
should ask about some form of  classical limit of the TMF description. 
What replaces the Chern chracter map should be a version for TMF 
(cf. the Miller character) of the Pontrjagin character map in KO-theory
\(
Ph: KO^*(X) \longrightarrow^{\!\!\!\!\!\!\!\!\!\!\otimes \C}~K^*(X)\longrightarrow^{\!\!\!\!\!\!\!\!\!\!\rm ch}
~~H^{**}(X;\Q)\;.
\)
The range is degree $4k$ cohomology, which indeed captures the 
field $G_4$ and its dual $G_8$ (or $\Theta$). We think of this as the 
combination 
$
{\sf v}~ G_4 + {\sf w}~ G_8
$,
for suitably identified generators ${\sf v}$ and ${\sf w}$, which make 
$H^{**}(X;\Q)$. Some aspects of this in connection to Spin K-theory
have been discussed in \cite{KSpin} (see also section \ref{int lift}).

\noindent {\bf 6.}
{\bf Infinite number of fields.}
We know that, physically, an element in differential K-theory corresponds 
to a collection of physical fields: all the RR fields. Analogously, an element in
(differential twisted) elliptic cohomology or (differential twisted)
$TMF$ should be given 
by a large collection of physical fields.  At the level of differential forms, the 
large number of fields 
should be seen in the same way that the K-theoretic RR fields are
seen at the level of differential forms via supergravity fields and their 
Hodge duals. 
But eleven-dimensional supergravity as traditionally known features only 
a single candidate  
field: the 3-form $C$-field (and its Hodge dual).
This certainly cannot model generic elements in $TMF$ by itself. Therefore, 
our previous discussion suggests a considerably richer structure hidden 
within and beyond eleven-dimensional supergravity. 
Recall that such a rich structure is also suggested by 
 hidden symmetries:

\vspace{2mm}
{\it  First,  at the level of differential forms:}

$\bullet$ In the dimensional reduction of eleven-dimensional 
supergravity (and hence type IIA supergravity) on tori 
$T^n$ with fluxes one gets the Cremmer-Julia \cite{CJ}
exceptional groups
$E_{n(n)}$, which are infinite-dimensional Kac-Moody 
groups for $n\geq 9$ (cf. \cite{J}). In the latter case, hence, there are an infinite
number of fields at the classical level. Passing to the quantum theory,
one has the $U$-duality groups $E_{n(n)}(\Z)$, the $\Z$-forms of the
above non-compact groups, and so we still have an infinite 
number of fields. 

$\bullet$ Already in eleven dimensions, the classic works of de Wit, Nicolai
( reviewed in \cite{dWN}),
and conjecture of \cite{Du}, imply the existence of Cremmer-Julia groups 
without compactification. One recent striking proposal is that of \cite{We}
in which the Lorentzian Kac-Moody algebra $E_{11}$ is proposed 
as a symmetry in M-theory. This proposal has withstood many checks.
The algebra $\frak{e}_{11}$ admits an infinite $\Z$-grading as
$
\ee_{11}= \cdots \oplus \ee_8 \oplus \cdots
$.

\vspace{2mm}
{\it Second, integrally:} The above fields should have refinements at the quantum level to 
whichever (generalized) cohomology theory ends up arising.

\noindent {\bf 7.} {\bf 
The topological term via higher classes.}
The topological part of the M-theory action is written in a suggestive compact 
form when lifted to the bounding twelve-manifold 
$Z^{12}$ in \cite{S1} \cite{S2} \cite{S3}. 
To do so, the {\it total 
String class} (in the notation of  those papers)
$\lambda=1+\lambda_1+\lambda_2+\cdots$ is introduced,
where $\lambda_1=p_1/2$ is the usual `$String$ class',
and $\lambda_2=p_2/2$ which is well-defined for Spin
manifolds. 
The interpretation of the class and the characters is
as degree $4k$ analogs
of the Chern class and the Chern character, mapping
from the cohomology theory describing M-theory to degree $4k$
cohomology. These are essentially the Spin characteristic classes
defined in \cite{Th} and were precursors to the discussion of 
Fivebrane structures in \cite{SSS1} \cite{SSS2} \cite{SSS3}.

\begin{theorem} 
{\cite{S1}} {\cite{S2}} {\cite{S3}}{\cite{S5}}
\begin{enumerate}
\item The M-theory fields are elements of a unified 
field strength.
The quantization formula on the total M-theory field reproduces 
the quantization
on $G_4$ and its dual $G_8$.
\item  $G_4$ can be viewed as an index 2-gerbe.
\end{enumerate}
\end{theorem}

\noindent {\bf 8. Twisted cohomology in M-theory.}
The characters can be extended to 
include the dual fields and to account for their dynamics. 
For the dual field, one
can pick either the straightforward degree seven Hodge
dual or its differential, the degree eight field $\Theta$. The
dual formulation favors the degree four/eight
combination whereas the duality-symmetric dynamics
favors the degree four/seven combination. 
The
second combination hints at a role for the prime $p=3$ in M-theory
analogous to the role played by the prime $p=2$
in K-theory, and consequently in string theory.
The combination leads to 
twisted (generalized) cohomology.
However, unlike the case in K-theory, the twist is given by 
a degree four class, suggesting relation to
TMF, since 
a twist of the latter
can be seen, at least heuristically,  
as a K-theory (degree three) twist on the loop space.

\begin{theorem} 
{\cite{S2}} {\cite{S3}} 
There is a twisted (graded) cohomology on the M-theory
fields with a twist given by a degree four class.
\end{theorem}


%


\noindent{\bf 9. The $q$-expansions.}
The $q$-expansions from the Miller character 
$
TMF \to K[[q]] \to H^*\Q [[q]]
$
will come from the comparison to type IIA string theory. 
We have the following diagram


\(
\xymatrix{
S^1
\ar[rr]
&&
\mathcal{W}^3~
\ar[d]^{\pi'}
\ar@{^{(}->}[rr]
&&
Y^{11}
\ar[d]^{\pi}
&&
S^1
\ar[ll]
 \ar@/_2pc/[llllll]_{=}
\\
&&
\Sigma_2~
\ar@{^{(}->}[rr]
&&
X^{10}
&&
}\;.
\)
There are two principal circle bundles: $\pi$ and $\pi'$. In the
dimensional reduction from M-theory to type IIA string theory
the two fibers are identified and from a membrane in eleven
dimensions we get a string in ten dimensions. 
At the level of 
partition functions of the targets (the fiber bundle $\pi$ in the diagram)
it was observed in \cite{KS1} that the resulting partition function,
formulated in elliptic cohomology, is a $q$-expansion of the 
K-theoretic partition function. There $q$ was built out of the generators
$v_1$ and $v_2$ at the prime 2, $q =v_1^3 v_2^{-1}$. This essentially
came from the fact that $EO_2({\rm pt})=\Z_2[[q]]$.

\vspace{3mm}
\noindent {\bf The $E$-theoretic quantization.}
In the untwisted case, the total field strength $F(x)$ of type II string theory as described
by K-theory 
is $2\pi$ times $\sqrt{\hat{A}(X)}ch(x)$.
The refinement of the description of the fields to elliptic cohomology is 
For any
elliptic cohomology theory $E$, there is a canonical map $E\to
K((q))$ (where $q$ is as above), so compose with the Chern
character to get a map $ch_E:E\to H^*((q))$.
The term $\sqrt{\hat{A}(X)}$ should be replaced by an
analogous term related to the Witten genus $\sigma(X)^{1/2}$ where 

\begin{theorem} [{\cite{KS3}}] 
{\it The $E$-theoretic quantization condition for the RR fields
is given by the formula for the elliptic field strength associated with $x$:
$F(x)=\sigma(X)^{1/2}ch_E(X)$, where $\sigma(X)$ is the characteristic class of $X$ associated with the power 
series $$\sigma(z)=(e^{z/2}-e^{-z/2}){\prod}_{n \geq 1}{}
\frac{(1-q^ne^z) (1-q^ne^{-z})}{(1-q^n)^2}\;.$$ }
\end{theorem}
The $\sigma$-function, in the $q\to 0$ limit, reduces to the
characteristic function of the $\hat{A}$-genus, thus reducing this
field strength to the type II field strength in the
$10$-dimensional limit.

\vspace{3mm}
The connection of this to charges is as follows. Considering the 
bundle $\pi'$, we get $q$-expansions on $\Sigma_2$ upon taking Fourier modes of the
circle. The boundary $\partial \Sigma_2$ has Chan-Paton charges on it.
Since the boundary of the string ends on a D-brane then certainly the
$q$-expansions will be seen by the D-brane as well. The result is that
the Chan-Paton bundle on the D-brane gets replaced by its
appropriate $q$-expansions 
$
E \mapsto S^{\bullet} E, \Lambda^{\bullet} E
$.
The charges of $D$-branes have been previously calculated at $q=1$,
i.e. at $\Lambda^{\bullet}=S^{\bullet}=1$,
as they are essentially the Chern character
\(
K^* (q=1) \longrightarrow H^* \Q (q=1)\;. 
\) 
The $q$-refinement of the D-brane charge formula should then be
\(
Q = \Phi_W (X^{10})ch_{ell}(E)\;,
\)
the Witten genus of $X^{10}$ twisted by the appropriate exterior and 
symmetric powers of the Chan-Paton bundle. 

\vspace{3mm}
\noindent{\bf The higher modes for the supermultipet.}
In an orthogonal
discussion to \cite{KS1},
we could also view $q$ as coming from 
Fourier modes on the circle, i.e. Kaluza-Klein modes and interpret $q$ accordingly.
Consider the supergravity fields $(g, C_3, \psi)$.
The coupling to vector bundles $V$ gives $V\otimes {\mathcal L}^k)$. 
For example $V$ is an $E_8$ bundle with characteristic
class $a$ and set $c_1({\mathcal L})$. Consider the connection 
$A$ on the vector bundle $V$. Coupling it to ${\mathcal L}^k$
leads to the connection $A \otimes e^{-ik\theta}$.
Then the C-field, which is essentially the Chern-Simons form
of the $E_8$ bundle, will also be also have Fourier components
as $C_3  \otimes e^{-ik \theta}$. This is of course compatible with,
and is in fact the same as, just the dimensional reduction of the
C-field directly from eleven dimensions without use of the $E_8$
bundle. 

\vspace{3mm}
\noindent {\bf 10. Concluding remark.}
 We have reviewed and indicated further development that a closer 
 examination of the deep structures involved in string- and M-theory
 indicates and shows that very rich cohomological phenomena are
 at work in the background. While it is now well-known that (differential)
 K-theory encodes much of the interesting structure of Type II string
 theory, we have shown and argued that by further reasoning along
 such lines one finds various generalized cohomological 
structures that go beyond K-theory. This includes Morava K-theory
and E-theory, elliptic cohomology, and $TMF$. In addition, we
have considered higher structures generalizing Spin structures,
such as String and Fivebrane structures, and have
emphasized a delicate interplay between all these.

\vspace{3mm}
{\bf \large Acknowledgements}

\vspace{2mm}
\noindent 
I would like to thank Matthew
Ando, David Lipsky, Corbett Redden, Urs Schreiber, and 
Jim Stasheff for  very useful 
discussions, and Dan Freed for useful remarks. I also thank the organizers 
of the CBMS  conference on ``Topology, $C^*$-algebras, and String Duality",
Robert Doran and Greg Friedman, for the kind invitation to the 
conference. 
Special thanks are due to the anonymous referee,
who kindly gave many helpful remarks and suggestions for
improving the presentation and strengthening the points
made in the article.  I also thank Arthur Greenspoon for 
editorial suggestions.


\end{document}